\documentclass[12pt,reqno]{amsart}

\newcommand\version{November 18, 2026}

\usepackage{amsmath,amsfonts,amsthm,amssymb,amsxtra}
\usepackage{bbm} % to get \1
\usepackage{hyperref} % creates hyperlinks to the references 

\newtheorem{theorem}{Theorem}%[section]
\newtheorem{proposition}[theorem]{Proposition}
\newtheorem{lemma}[theorem]{Lemma}

\theoremstyle{definition}

\theoremstyle{remark}

\newtheorem{remark}[theorem]{Remark}
\newtheorem{remarks}[theorem]{Remarks}

\newcommand{\1}{\mathbbm{1}}

\newcommand{\C}{\mathbb{C}}

\renewcommand{\epsilon}{\varepsilon}

\newcommand{\F}{\mathcal{F}}
\newcommand{\ii}{\mathrm i}
\newcommand{\loc}{{\rm loc}}

\renewcommand{\phi}{\varphi}
\newcommand{\R}{\mathbb{R}}

\newcommand{\Sp}{\mathcal S}
\newcommand{\Sch}{\mathcal T}
\newcommand{\T}{\mathbb{T}}

\newcommand{\Z}{\mathbb{Z}}

\newcommand{\dd}{\,\mathrm d}

\DeclareMathOperator{\supp}{supp}
\DeclareMathOperator{\sgn}{sgn}
\DeclareMathOperator{\Tr}{Tr}

\begin{document}

\title[Cwikel and Kato--Seiler--Simon estimates for commutators]{Cwikel and Kato--Seiler--Simon estimates\\ for commutators}

\author{Rupert L. Frank}
\address[Rupert L. Frank]{Mathe\-matisches Institut, Ludwig-Maximilians Universit\"at M\"unchen, The\-resienstr.~39, 80333 M\"un\-chen, Germany, and Munich Center for Quantum Science and Technology, Schel\-ling\-str.~4, 80799 M\"unchen, Germany.}
\email{r.frank@lmu.de}

\thanks{\copyright\, 2026 by the author. This paper may be reproduced, in its entirety, for non-commercial purposes.\\
	The author acknowledges partial support from the German Research Foundation through EXC-2111-390814868 and TRR 352-Project-ID 470903074 and from OpenAI through the ChatGPT for Academic Researchers program.}

\begin{abstract}
	We prove that the commutator $[a(X),b(-\ii\nabla)]$ belongs to the weak Schatten class $\Sp^p_{\rm weak}$ provided $\nabla a\in L^p$ and $\nabla b\in L^p_{\rm weak}$ for $\max\{1,2d/(d+2)\}<p<\infty$. This is an analogue of Cwikel's estimate, with a larger range of exponents. We also prove analogues of the Kato--Seiler--Simon and Birman--Solomyak--Simon estimates. Counterexamples show the necessity of the lower thresholds and the failure of the Lipschitz endpoint estimate.
\end{abstract}

\date{\version}

\maketitle

%%%%%%%%%%%%%%%%%%%%%%%%%%%%%%%%
%%%%%%%%%%%%%%%%%%%%%%%%%%%%%%%%

\section{Introduction and main result}

A fundamental result in spectral theory, which was conjectured by Simon \cite{Simon76} and proved by Cwikel \cite{Cwikel}, states that,
for $d\geq1$ and $2<p<\infty$, the assumptions $a\in L^p(\R^d)$ and
$b\in L^p_{\rm weak}(\R^d)$ imply $a(X)b(-\ii\nabla)\in\Sp^p_{\rm weak}$, with
\begin{equation}\label{eq:cwikel}
\| a(X) b(-\ii\nabla) \|_{\Sp^p_{\rm weak}} \lesssim_{d,p} \|a\|_{L^p} \|b\|_{L^p_{\rm weak}}.
\end{equation}
A precursor of Cwikel's estimate is due to Kato (unpublished) and
Seiler--Simon \cite{SeSi} and states that, for $d\geq1$ and
$2\leq p<\infty$, the assumptions $a,b\in L^p(\R^d)$ imply
$a(X)b(-\ii\nabla)\in\Sp^p$, with
\begin{equation}\label{eq:kss}
	\| a(X) b(-\ii\nabla) \|_{\Sp^p} \lesssim_{d,p} \|a\|_{L^p} \|b\|_{L^p}.
\end{equation}

The notation \(\lesssim_{d,p}\) means that the inequality holds up to a multiplicative constant depending only on \(d\) and \(p\). Moreover, $L^p(\R^d)$ denotes the usual $L^p$-space on $\R^d$ and $L^p_{\rm weak}(\R^d)$ denotes the corresponding weak $L^p$-space, consisting of all measurable functions $f$ on $\R^d$ for which the quasinorm
\[
\|f\|_{L^p_{\rm weak}}:=\sup_{t>0}t\,|\{|f|>t\}|^{1/p}
\]
is finite. Similarly, $\Sp^p=\Sp^p(L^2(\R^d))$ denotes the Schatten class, consisting of all compact operators $T$ in $L^2(\R^d)$ whose singular values $s_n(T)$, $n\geq 1$, belong to $\ell^p$ with its natural quasinorm (a norm if $p\geq 1$), and $\Sp^p_{\rm weak}=\Sp^p_{\rm weak}(L^2(\R^d))$ denotes the weak Schatten class, consisting of all compact operators $T$ in $L^2(\R^d)$ for which the quasinorm
\[
\|T\|_{\Sp^p_{\rm weak}}:=\sup_{n\geq1}n^{1/p}s_n(T)
\]
is finite.

The operator $a(X)$ acts as multiplication by the function $a$ in position space, while $b(-\ii\nabla)$ acts as multiplication by the function $b$ in momentum space, that is, $\mathcal F b(-\ii\nabla)\mathcal F^*$ is the operator of multiplication by $b$, where $\mathcal F$ denotes Fourier transform. For a precise definition of $a(X) b(-\ii\nabla)$ as a distributional form, see Section~\ref{sec:commid}.

As we will discuss later in this introduction, both the inequalities of Cwikel and of Kato--Seiler--Simon capture a semiclassical intuition and therefore play a crucial role in the study of spectral asymptotics of differential and pseudo-differential operators. Cwikel's estimate implies, in particular, the Cwikel--Lieb--Rozenblum inequality \cite{Cwikel,Li,Ro} for the number of negative eigenvalues of Schr\"odinger operators and allows one to prove semiclassical spectral asymptotics under minimal regularity assumptions. The Kato--Seiler--Simon inequality can be used to derive non-endpoint Lieb--Thirring inequalities \cite{LT}. We refer to \cite{Fr2,FrLaWe} for references and background on this.

Our goal here is to prove an analogue of these inequalities for commutators. The following is our main result.

\begin{theorem}\label{thm:main}
	Let $d\geq1$ and $\max\{1,2d/(d+2)\}<p<\infty$.
	\begin{enumerate}
		\item[(a)] If $a,b\in L^1_{\mathrm{loc}}(\R^d)$
		have distributional gradients satisfying $\nabla a\in L^p(\R^d)$, $\nabla b\in L^p_{\rm weak}(\R^d)$, then $[a(X),b(-\ii\nabla)]\in \Sp^p_{\rm weak}$ and
		\begin{equation}\label{eq:main}
			\|[a(X),b(-\ii\nabla)]\|_{\Sp^p_{\rm weak}}
			\lesssim_{d,p} \|\nabla a\|_{L^p}\|\nabla b\|_{L^p_{\rm weak}}.
		\end{equation}
		
		\item[(b)] If $a,b\in L^1_{\mathrm{loc}}(\R^d)$
		have distributional gradients satisfying $\nabla a,\nabla b\in L^p(\R^d)$, then $[a(X),b(-\ii\nabla)]\in\Sp^p$ and 
		\begin{equation}\label{eq:strong}
			\|[a(X),b(-\ii\nabla)]\|_{\Sp^p}
			\lesssim_{d,p} \|\nabla a\|_{L^p}\|\nabla b\|_{L^p}.
		\end{equation}
	\end{enumerate}
\end{theorem}

\begin{remarks}\label{rem:main}
	(a) More precisely, we show that the commutator, initially interpreted as a distributional form, extends to an operator in the appropriate Schatten class. This is further discussed in Section \ref{sec:commid}.\\
	(b) In the special case $b(\xi) = |\xi|^{1-d/p}$ with $p\ne d$, part (a) was proved in our earlier work \cite{FSZ} for the same range of \(p\). The new ingredients needed for general $b$ are explained at the end of the introduction.\\
	(c) Part (b) in the special case $p=2$ is due to Simon \cite{Simon79}. In this case, one obtains the bound with explicit constant
	\begin{equation}
		\label{eq:strong2}
		\|[a(X),b(-\ii\nabla)]\|_{\Sp^2}
		\leq (2\pi)^{-d/2} \|\nabla a\|_{L^2}\|\nabla b\|_{L^2};
	\end{equation}
	see Remark \ref{rem:simon}. (Simon’s displayed proof appears to contain a normalization slip.)\\
	(d) For $d\leq2$, the lower restriction $p>1$ in (a) and (b) is optimal: Theorem~\ref{thm:counter} shows that even weaker estimates fail for every $0<p\leq1$.\\
	(e) For $d\geq 3$, the assumption $p>2d/(d+2)$ in (a) is best possible. Indeed, we take a nonzero real $a\in C_c^\infty(\R^d)$ and $b(\xi):=|\xi|^{-d/2}$. Then $\nabla b\in L^{2d/(d+2)}_{\rm weak}(\R^d)$, but the commutator, defined distributionally on $C_c^\infty(\R^d)$, has no bounded extension to $L^2(\R^d)$, by \cite[Lemma~3.3]{FSZ}.
	To recall the reason, we choose $\psi\in C_c^\infty(\R^d)$ with $\int a\psi\,dx\ne0$.
	Outside the support of $a$, the commutator applied to $\psi$, interpreted as a distribution, has a nonzero leading term proportional to $|x|^{-d/2}$ and hence is not in $L^2(\R^d)$.\\
	(f) For $d\geq3$, the question whether part (b) extends to $p=2d/(d+2)$ remains open. It is easy to see, however, that it cannot be extended to $1<p<2d/(d+2)$. Indeed, choosing $\chi\in C_c^\infty(\R^d)$ equal to one near the origin, we see that $a(x):=\chi(x)|x|^{-d/2}$ has distributional gradient in $L^p$, but
	$a\notin L^2_{\loc}$. We choose a nonzero real $b\in C_c^\infty$ and a
	test function $\psi$ supported away from the origin such that
	$b(-\ii\nabla)\psi(0)\ne0$. The distribution
	$a\,b(-\ii\nabla)\psi$ is not locally in $L^2$, whereas
	$b(-\ii\nabla)(a\psi)$ is smooth. Thus this distributional commutator has no bounded extension to $L^2(\R^d)$,
	even though both gradients belong to $L^p$ for any $1<p<2d/(d+2)$. The cases $p\leq1$ are covered by Theorem~\ref{thm:counter}.\\
	(g) The assumption $p<\infty$ is best possible in the sense that $\nabla a,\nabla b\in L^\infty(\R^d)$ does not guarantee that $[a(X),b(-\ii\nabla)]$ is bounded; see Theorem \ref{thm:counterexample} and Remark~\ref{rem:unbounded} below.
\end{remarks}

%%%%%%%%%%%%%%%%%%%%

\subsection*{Estimates for smaller exponents}

The Cwikel and Kato--Seiler--Simon inequalities require the assumptions $p>2$ and $p\geq 2$, respectively. Even earlier, Birman and Solomyak \cite{BS} had obtained bounds for $p\leq 2$. To state these, we recall the definition of the spaces they introduced. Letting $Q_k:=k+[0,1)^d$ for all $k\in\mathbb Z^d$, for a measurable scalar or vector-valued function $f$ on $\R^d$, we set
\[
f_k:=\left(\int_{Q_k}|f(x)|^2\,dx\right)^{1/2},\qquad
\|f\|_{\ell^p(L^2)}:=\|(f_k)\|_{\ell^p},\qquad
\|f\|_{\ell^p_{\rm weak}(L^2)}:=\|(f_k)\|_{\ell^p_{\rm weak}},
\]
where the weak sequence norm is
$\|(c_k)\|_{\ell^p_{\rm weak}} :=\sup_{t>0}t\,(\#\{k:|c_k|>t\})^{1/p}$.
We note that
\[
\begin{aligned}
\|f\|_{L^p}&\leq \|f\|_{\ell^p(L^2)}
&&\text{for }0<p\leq2,\\
\|f\|_{L^p_{\rm weak}}&\lesssim_p \|f\|_{\ell^p_{\rm weak}(L^2)}
&&\text{for }0<p<2.
\end{aligned}
\]
The first bound follows from H\"older's inequality on each unit cube. For the second, Chebyshev's inequality gives, for every $t>0$,
\[
|\{|f|>t\}|\leq\sum_k\min\{1,t^{-2}f_k^2\}
=2t^{-2}\int_0^t s\,\#\{k:f_k>s\}\,ds
\leq\frac{2}{2-p}\,t^{-p}\|f\|_{\ell^p_{\rm weak}(L^2)}^p.
\]

The Birman--Solomyak bound \cite{BS} states that, for $d\geq1$ and
$0<p\leq2$, the assumptions $a,b\in\ell^p(L^2)$ imply
$a(X)b(-\ii\nabla)\in\Sp^p$, with
\begin{equation}\label{eq:bs}
	\| a(X) b(-\ii\nabla) \|_{\Sp^p} \lesssim_{d,p} \|a\|_{\ell^p(L^2)} \|b\|_{\ell^p(L^2)} \,.
\end{equation}
The weak version of this estimate is due to Simon \cite[Theorem~4.6]{Simon} (for $1<p<2$) and Birman--Karadzhov--Solomyak \cite[Section~5.7]{BKS} (for $0<p<2$). It states that, for $d\geq1$ and $0<p<2$, the assumptions
$a\in\ell^p(L^2)$ and $b\in\ell^p_{\rm weak}(L^2)$ imply
$a(X)b(-\ii\nabla)\in\Sp^p_{\rm weak}$, with
\begin{equation}\label{eq:bsweak}
	\| a(X) b(-\ii\nabla) \|_{\Sp^p_{\rm weak}} \lesssim_{d,p} \|a\|_{\ell^p(L^2)} \|b\|_{\ell^p_{\rm weak}(L^2)} \,.
\end{equation}

Our next theorem gives commutator bounds when the gradients belong to Birman--Solomyak spaces. These bounds add to Theorem~\ref{thm:main} only in the range $1<p\leq2d/(d+2)$, which is nonempty only for $d\geq3$; for larger $p$, they follow from that theorem and the embeddings above.

\begin{theorem}\label{thm:bscomm}
	Let $d\geq1$ and $1<p<2$.
	\begin{enumerate}
		\item[(a)] If $a,b\in L^1_{\mathrm{loc}}(\R^d)$
		have distributional gradients satisfying $\nabla a\in \ell^p(L^2)$, $\nabla b\in \ell^p_{\rm weak}(L^2)$, then $[a(X),b(-\ii\nabla)]\in \Sp^p_{\rm weak}$ and
		\begin{equation}\label{eq:bscomm1}
			\|[a(X),b(-\ii\nabla)]\|_{\Sp^p_{\rm weak}}
			\lesssim_{d,p} \|\nabla a\|_{\ell^p(L^2)}\|\nabla b\|_{\ell^p_{\rm weak}(L^2)}.
		\end{equation}
		
		\item[(b)] If $a,b\in L^1_{\mathrm{loc}}(\R^d)$
		have distributional gradients satisfying $\nabla a, \nabla b\in \ell^p(L^2)$, then $[a(X),b(-\ii\nabla)]\in \Sp^p$ and
		\begin{equation}\label{eq:bscomm2}
			\|[a(X),b(-\ii\nabla)]\|_{\Sp^p}
			\lesssim_{d,p} \|\nabla a\|_{\ell^p(L^2)}\|\nabla b\|_{\ell^p(L^2)}.
		\end{equation}
	\end{enumerate}
\end{theorem}

\begin{remarks}\label{rem:bscomm}
	(a) As in Theorem \ref{thm:main}, the commutator is initially interpreted as a distributional form; see Section \ref{sec:commid}.\\
	(b) The assumption $p>1$ is best possible, as shown in Theorem~\ref{thm:counter}.
\end{remarks}

%%%%%%%%%%%%%%%%%%%%%

\subsection*{Counterexamples}

The next two results show that the restrictions $p<\infty$ and $p>1$ cannot be relaxed.

The restriction $p<\infty$ in Theorem~\ref{thm:main} cannot be removed. This gives a negative answer to a question raised by Davies \cite[Footnote 12 on Problem 5.7.2]{Davies}.

\begin{theorem}\label{thm:counterexample}
	For every $d\geq1$, there is no implied constant such that the estimate
	\[
	\|[a(X),b(-\ii\nabla)]\|
	\lesssim_d\|\nabla a\|_\infty\|\nabla b\|_\infty
	\]
	holds for all real $a,b\in C_c^\infty(\R^d)$.
\end{theorem}

As we have recalled, for product operators $a(X)b(-\ii\nabla)$, there are also Schatten class bounds for $0<p\leq1$. The following result shows that the corresponding commutator estimates fail in every dimension.

\begin{theorem}\label{thm:counter}
	For every $d\geq1$ and $0<p\leq1$, neither of the estimates
	\begin{align*}
	\|[a(X),b(-\ii\nabla)]\|_{\Sp^p}
	&\lesssim_{d,p}\|\nabla a\|_{\ell^p(L^2)}
	                 \|\nabla b\|_{\ell^p(L^2)},\\
	\|[a(X),b(-\ii\nabla)]\|_{\Sp^p_{\rm weak}}
	&\lesssim_{d,p}\|\nabla a\|_{\ell^p(L^2)}
	                 \|\nabla b\|_{\ell^p_{\rm weak}(L^2)}
	\end{align*}
	holds with an implied constant independent of the real functions $a,b\in C_c^\infty(\R^d)$.
\end{theorem}

%%%%%%%%%%%%%%%%%%%%%

\subsection*{Extensions}

In this paper we have tried to strike a balance between optimal results and reasonably simple proofs. There are several possible extensions that we have not pursued, but which could be interesting and useful in applications. Those include
\begin{itemize}
	\item analogues of the $\Sp^2_{\rm weak}$-bound from \cite{LSZ},
	\item bounds in general Lorentz--Schatten ideals $\Sp^{p,q}$ as in \cite{BKS},
	\item bounds in non-power ideals as in \cite{We},
	\item bounds in an operator-valued setting as in \cite{Hu},
	\item bounds beyond the Euclidean setting, as in \cite{BKS}.
\end{itemize}
We have also not paid attention to the optimal (or even good) constants in our commutator bounds. The optimal constant in the Kato--Seiler--Simon inequality is the semiclassical constant. Rather tight bounds on the constant in Cwikel's inequality were found in \cite{Fr,HKRV}.

%%%%%%%%%%%%%%%%%%%%%

\subsection*{The semiclassical intuition}

Let us provide some background on why the Cwikel estimate and its new analogue for commutators capture a semiclassical structure. For this heuristic discussion, we assume that $a$ and $b$ are smooth and satisfy the usual assumptions for the computation of principal symbols. We interpret the phase-space expressions formally when $a$ or $b$ is singular.

We have
\[
\begin{aligned}
\| a(X)b(-\ii\nabla) \|_{\Sp^p_{\rm weak}}^p
&= \sup_{\lambda>0} \lambda^p n(\lambda,a(X)b(-\ii\nabla))\\
&= \sup_{\kappa>0} \kappa^{p/2} n(\kappa,a(X)|b(-\ii\nabla)|^2 \overline{a}(X)),
\end{aligned}
\]
where $n(\lambda,T)$ is the number of singular values of a compact operator $T$ that exceed $\lambda$. The operator $a(X)|b(-\ii\nabla)|^2 \overline{a}(X)$ has principal symbol $|a(x)|^2 |b(\xi)|^2$, so a (non-rigorous) semiclassical approximation according to Planck's principle (or Weyl's law) suggests
\[
n(\kappa,a(X)|b(-\ii\nabla)|^2 \overline{a}(X)) \approx \iint_{\R^d\times\R^d} \1(|a(x)|^2 |b(\xi)|^2>\kappa) \,\frac{dx\,d\xi}{(2\pi)^d}.
\]
We can use the definition of $\|b\|_{L^p_{\rm weak}}$ to bound the double integral
\begin{align*}
	\iint_{\R^d\times\R^d} \1(|a(x)|^2 |b(\xi)|^2>\kappa) \,\frac{dx\,d\xi}{(2\pi)^d}
	& \leq (2\pi)^{-d} \kappa^{-p/2} \|b\|_{L^p_{\rm weak}}^p \int_{\R^d} |a(x)|^p\,dx \\
	& = (2\pi)^{-d} \kappa^{-p/2} \|a\|_{L^p}^p \|b\|_{L^p_{\rm weak}}^p.
\end{align*}
While this is only an upper bound, it is a rather sharp one, since it turns into an equality in the special case $b(\xi)=|\xi|^{-d/p}$. In this sense, Cwikel's estimate turns a non-rigorous semiclassical approximation into a rigorous upper bound.

The argument for the Kato--Seiler--Simon inequality is similar, except now we approximate
\[
\begin{aligned}
\|a(X)b(-\ii\nabla)\|_{\Sp^p}^p
&=\Tr \left( a(X)|b(-\ii\nabla)|^2 \overline{a}(X) \right)^{p/2}\\
&\approx\iint_{\R^d\times\R^d}|a(x)|^p|b(\xi)|^p\,\frac{dx\,d\xi}{(2\pi)^d}.
\end{aligned}
\]

Let us turn to commutators. The operator $[a(X),b(-\ii\nabla)]$ has principal symbol $\ii \nabla a(x)\cdot\nabla b(\xi)$, so the square of its absolute value has principal symbol $|\nabla a(x)\cdot \nabla b(\xi)|^2$. Thus, in the expression
\[
\| [a(X),b(-\ii\nabla)] \|_{\Sp^p_{\rm weak}}^p = \sup_{\kappa>0} \kappa^{p/2} n(\kappa,|[a(X),b(-\ii\nabla)]|^2)
\]
we may approximate
\[
n(\kappa,|[a(X),b(-\ii\nabla)]|^2) \approx \iint_{\R^d\times\R^d} \1(|\nabla a(x)\cdot\nabla b(\xi)|^2 >\kappa) \,\frac{dx\,d\xi}{(2\pi)^d}.
\]
The bound $|\nabla a(x)\cdot\nabla b(\xi)|\leq |\nabla a(x)||\nabla b(\xi)|$ and the same argument as before give
\begin{align*}
	\iint_{\R^d\times\R^d} \1(|\nabla a(x)\cdot\nabla b(\xi)|^2 >\kappa) \,\frac{dx\,d\xi}{(2\pi)^d}
	& \leq (2\pi)^{-d} \kappa^{-p/2} \|\nabla a\|_{L^p}^p \|\nabla b\|_{L^p_{\rm weak}}^p.
\end{align*}
This is optimal up to a multiplicative constant depending on $d$ and $p$ when $b(\xi) = |\xi|^{1-d/p}$ and $p\ne d$. Thus, we see that our new Cwikel-type bound captures the semiclassical approximation for commutators. The same applies to our new Kato--Seiler--Simon bound.

The above heuristics already suggest that there is a relation between our new commutator estimates for $[a(X),b(-\ii\nabla)]$ and existing estimates for $|\nabla a|(X)|\nabla b|(-\ii\nabla)$. This is what our proof makes rigorous.

%%%%%%%%%%%%%%%%%%%%%

\subsection*{Outline of the paper}

Let us briefly describe the strategy in this paper. The positive results, Theorems \ref{thm:main} and \ref{thm:bscomm}, are proved in Section \ref{sec:commid}. For Theorem \ref{thm:main}, the fundamental tool is identity \eqref{eq:weighted-identity}, which represents the commutator $[a(X),b(-\ii\nabla)]$ as a sum of terms obtained by applying Schur multipliers to products of the form
\[
g^\theta(-\ii\nabla)(\partial_j a)(X)g^{1-\theta}(-\ii\nabla)
\]
for a suitably chosen function $g$ and a parameter $0\leq\theta\leq1$. In fact, a refinement of the Haj\l asz inequality allows us to choose $g$ and $\theta$ so that $g$ is controlled by $|\nabla b|$ in the relevant norm and both $p/\theta$ and $p/(1-\theta)$ exceed two. Factoring each such term into two products allows us to apply Cwikel's inequality, or the Kato--Seiler--Simon inequality, at these larger exponents. For the Birman--Solomyak bounds we take $\theta=1$ in the same identity.

An essential step in this argument is the proof that the relevant Schur multipliers are bounded on the Schatten space of interest. This is deduced in Section \ref{sec:multiplier} from \cite{CGPT}.

The same combination of Schur multipliers and factorizations underlies our previous work \cite{FSZ}. There, however, the power structure of $b$ leads to a different representation. The refined Haj\l asz inequality in Lemma~\ref{lem:hajlasz-gradients} permits us to treat arbitrary $b$ with $\nabla b\in L^p_{\rm weak}$, while retaining the same range of $p$ as in \cite{FSZ}.

The remaining sections deal with the counterexamples. In Section~\ref{sec:endpoint}, to disprove a bound for $p=\infty$ (Theorem~\ref{thm:counterexample}), we choose oscillations in $a$ and $b$ whose derivatives stay bounded through cancellation, while the terms that survive in the commutator produce an arbitrarily large norm. In Section~\ref{sec:counter}, to disprove bounds for $p\leq 1$ (Theorem~\ref{thm:counter}), we use oscillations at well separated frequencies whose contributions to the trace norm accumulate while the relevant gradient norms remain controlled. The weak estimate at $p=1$ requires a different choice of $b$. In each case, we first construct examples in one dimension and then extend them to arbitrary dimension.

%%%%%%%%%%%%%%%%%%%%%%%%%

\section{Boundedness of Schur multipliers}\label{sec:multiplier}

For a bounded measurable function $m$ on $\R^d\times\R^d$, we write $\Sch_m$ for the corresponding Schur multiplier, that is, the operator of multiplying the integral kernel of an operator on $L^2(\R^d)$ by $m$. This is a bounded operator on the Hilbert--Schmidt class $\Sp^2$. The following result, which is the main result of this section, says that for certain $m$ the multiplier $\Sch_m$ is bounded on all Schatten spaces $\Sp^p$, $1<p<\infty$, and their weak analogues. Whenever $\Sch_m$ is bounded on $\Sp^{p_0}$ and $\Sp^{p_1}$ with $1<p_0<p<p_1<\infty$, its action on $\Sp^p_{\rm weak}$ means the restriction of its compatible action on $\Sp^{p_0}+\Sp^{p_1}$.

\begin{proposition}\label{prop:weighted-multiplier}
We assume that $0\leq\theta\leq1$, that
$b:\R^d\to\C$ and $g:\R^d\to[0,\infty)$ are measurable and finite
almost everywhere, and that
\begin{equation}\label{eq:weighted-hajlasz}
|b(\xi)-b(\eta)|\leq |\xi-\eta|
\bigl(g(\xi)^\theta g(\eta)^{1-\theta}+g(\xi)^{1-\theta} g(\eta)^\theta\bigr)
\qquad\text{for a.e.}\ (\xi,\eta)\in\R^d\times\R^d.
\end{equation}
We define
\[
r_\theta(\xi,\eta)
:=\frac{b(\xi)-b(\eta)}
{|\xi-\eta|\bigl(g(\xi)^\theta g(\eta)^{1-\theta}
                 +g(\xi)^{1-\theta} g(\eta)^\theta\bigr)},
\]
with value zero when the denominator vanishes. Then, for $1<p<\infty$,
\[
\|\Sch_{r_\theta}\|_{\Sp^p\to\Sp^p}
\lesssim_{d,p}1,\qquad
\|\Sch_{r_\theta}\|_{\Sp^p_{\rm weak}\to\Sp^p_{\rm weak}}
\lesssim_{d,p}1.
\]
\end{proposition}

We emphasize that the implied constants in the conclusion are independent of $b$, $g$ and $\theta$.

In \eqref{eq:weighted-hajlasz} and everywhere below, we interpret
$g^0\equiv1$, including on the set where $g=0$. 

For $\theta=1$, the hypothesis \eqref{eq:weighted-hajlasz} becomes
\begin{equation}\label{eq:hajlasz}
|b(\xi)-b(\eta)|\leq|\xi-\eta|\bigl(g(\xi)+g(\eta)\bigr)
\qquad\text{for a.e.}\ (\xi,\eta)\in\R^d\times\R^d.
\end{equation}
A function $g$ satisfying \eqref{eq:hajlasz} for a given $b$ is called
a Haj\l asz gradient of $b$. This is an important concept in the theory
of Sobolev spaces on metric measure spaces, but here we will not need
anything from this theory.

The proof of this proposition relies on a H\"ormander--Mikhlin theorem for Schur multipliers of difference form $m(u,v)=\phi(u-v)$. We restate this result for easier reference following \cite[Introduction (HMS)]{CGPT}. For $\phi\in C^{N}(\R^n\setminus\{0\})$, we define
\begin{equation}\label{eq:mikhlin}
	\mathfrak M_{N}(\phi):=
	\sum_{|\alpha|\leq N}\sup_{z\ne0}
	|z|^{|\alpha|}|\partial^\alpha\phi(z)| \,.
\end{equation}

\begin{theorem}
	\label{thm:hms}
	For $n\geq1$ and $1<p<\infty$, if $\phi\in C^{\lfloor n/2\rfloor+1}(\R^n\setminus\{0\})$,
	then
	\begin{equation}\label{eq:hms}
		\|\Sch_{\phi(u-v)}\|_{\Sp^p(L^2(\R^n))\to\Sp^p(L^2(\R^n))}
		\lesssim_{n,p} \mathfrak M_{\lfloor n/2\rfloor+1}(\phi).
	\end{equation}
\end{theorem}

In fact, \cite{CGPT} gives a stronger, completely bounded estimate, but the stated
ordinary Schatten bound suffices for our purposes. 

To derive Proposition \ref{prop:weighted-multiplier} from Theorem \ref{thm:hms}, we need the following result about Schur multipliers.

\begin{lemma}\label{lem:pullback}
	For $1<p<\infty$, we assume that $M_k$ are bounded continuous functions on
	$\R^n\times\R^n$ with
	\[
	C:=\sup_k\|\Sch_{M_k}\|_{\Sp^p(L^2(\R^n))\to\Sp^p(L^2(\R^n))} <\infty 
	\]
	and that $M_k(u,v)\to M(u,v)$ for every $(u,v)\in\R^n\times\R^n$. Then for any measurable maps $\alpha,\beta:\R^d\to\R^n$, the function
	\[
	(\xi,\eta)\longmapsto M(\alpha(\xi),\beta(\eta))
	\]
	defines a Schur multiplier on $\Sp^p(L^2(\R^d))$ with norm at most $C$.
\end{lemma}

\begin{proof}
	We first consider a single bounded continuous kernel $M$ with multiplier
	norm at most $C$. We choose finitely many points
	$\lambda_1,\ldots,\lambda_r$ and $\mu_1,\ldots,\mu_t$
	in $\R^n$, allowing repetitions, and orthonormal families $(e_i^\delta)_{i=1}^r$ and
	$(f_j^\delta)_{j=1}^t$, with $e_i^\delta$ supported in
	$B_\delta(\lambda_i)$ and $f_j^\delta$ supported in
	$B_\delta(\mu_j)$. For distinct points the
	balls are disjoint when $\delta$ is small; for repeated points we choose
	several orthonormal functions in the same ball. We stress that the two families need not be orthogonal to each other.
	We denote by $J_\delta:\mathbb C^r\to L^2(\R^n)$ and
	$H_\delta:\mathbb C^t\to L^2(\R^n)$ the corresponding isometries.
	For every matrix $A=(a_{ij})\in\C^{r\times t}$, continuity of $M$ gives
	\[
	J_\delta^*\Sch_M(J_\delta A H_\delta^*)H_\delta
	\longrightarrow \bigl(M(\lambda_i,\mu_j)a_{ij}\bigr)_{i,j}
	\qquad\text{as }\delta\downarrow0.
	\]
	Indeed, we replace $M(u,v)$ on the support of each
	$e_i^\delta(u)\overline{f_j^\delta(v)}$ by $M(\lambda_i,\mu_j)$.
	The resulting error tends to zero in Hilbert--Schmidt norm, since
	there are only finitely many pairs. Embedding by isometries preserves
	the nonzero singular values, and multiplication by their adjoints
	does not increase a Schatten norm. Hence
	\begin{equation}\label{eq:finitepull}
		\|\bigl(M(\lambda_i,\mu_j)a_{ij}\bigr)\|_{\Sp^p}
		\leq C\|A\|_{\Sp^p}.
	\end{equation}
	Taking a one-entry matrix gives $|M(u,v)|\leq C$ everywhere.
	Repeating the points in \eqref{eq:finitepull} permits matrix blocks of any
	finite size.
	
	For finite-valued measurable maps $\alpha,\beta$, we write $P_i,Q_j$ for the
	multiplication projections onto the level sets with values $\lambda_i,\mu_j$,
	respectively. The corresponding multiplier is
	\[
	T\longmapsto\sum_{i,j}M(\lambda_i,\mu_j)P_iTQ_j.
	\]
	For finite-dimensional projection ranges its norm is at most $C$ by
	\eqref{eq:finitepull}. We choose finite-rank projections
	$P_{i,\nu}\leq P_i$ and $Q_{j,\nu}\leq Q_j$ converging strongly to
	$P_i$ and $Q_j$, respectively. For every $T\in\Sp^p$, we have
	\[
	P_{i,\nu}TQ_{j,\nu}\longrightarrow P_iTQ_j\quad\text{in }\Sp^p.
	\]
	Thus the same estimate holds
	for the level-set subspaces of $L^2(\R^d)$, even when these
	subspaces are infinite dimensional.
	
	For general $\alpha,\beta$, we choose finite-valued measurable maps $\alpha_\ell$, $\beta_\ell$ such that $\alpha_\ell\to\alpha$ and $\beta_\ell\to\beta$ pointwise. If $T$ is finite rank with Hilbert--Schmidt kernel $A$,
	then dominated convergence, using $|M|\leq C$, gives
	\[
	M(\alpha_\ell(\xi),\beta_\ell(\eta))A(\xi,\eta)
	\longrightarrow M(\alpha(\xi),\beta(\eta))A(\xi,\eta)
	\quad\hbox{in } L^2(\R^d\times\R^d).
	\]
	The corresponding operators converge in $\Sp^2$, and their
	$\Sp^p$ norms are at most $C\|T\|_{\Sp^p}$.
	Convergence in $\Sp^2$ implies convergence in operator norm, so the
	non-commutative Fatou lemma \cite[Theorem~2.7(d)]{Simon} gives the
	same bound for the limit. We extend by density from finite-rank operators.
	
	Finally, we apply this argument to each $M_k$. The pointwise bounds
	$|M_k|\leq C$ give dominated convergence in $\Sp^2$ after substitution;
	the same Fatou and density argument proves the claim
	for $M$.
\end{proof}

Next, we provide some examples of Schur multipliers. This is where Theorem \ref{thm:hms} enters. For $s\in\R$ we put
\[
\Phi_{n,s}(z):=
\begin{cases}|z|^{is},&z\ne0,\\0,&z=0,\end{cases}
\qquad
\Psi_s(t):=
\begin{cases}\sgn(t)|t|^{is},&t\ne0,\\0,&t=0.\end{cases}
\]

\begin{lemma}
	\label{lem:powers}
	For $1<p<\infty$, $s\in\R$, and measurable, almost everywhere finite
	$b:\R^d\to\R$ and $g:\R^d\to[0,\infty)$, we have
	\begin{align}
		\|\Sch_{\Phi_{d,s}(\xi-\eta)}\|_{\Sp^p\to\Sp^p}
		&\lesssim_{d,p}(1+|s|)^{\lfloor d/2\rfloor+1},\label{eq:power-euclidean}\\
		\|\Sch_{\Psi_s(b(\xi)-b(\eta))}\|_{\Sp^p\to\Sp^p}
		&\lesssim_p1+|s|,\label{eq:power-signed}\\
		\|\Sch_{\Phi_{1,s}(g(\xi)+g(\eta))}\|_{\Sp^p\to\Sp^p}
		&\lesssim_p1+|s|.\label{eq:power-sum}
	\end{align}
\end{lemma}

\begin{proof}
	For a multi-index $\alpha$, differentiation gives
	\[
	\partial^\alpha |z|^{is}
	=|z|^{is-|\alpha|}P_{\alpha,s}(z/|z|),\qquad
	\sup_{\omega\in S^{n-1}}|P_{\alpha,s}(\omega)|
	\lesssim_\alpha(1+|s|)^{|\alpha|}.
	\]
	Indeed, induction gives a polynomial in $s$ of degree at most $|\alpha|$,
	with smooth coefficients on the sphere. Thus
	$\mathfrak M_{\lfloor n/2\rfloor+1}(|z|^{is})\lesssim_n(1+|s|)^{\lfloor n/2\rfloor+1}$, proving
	\eqref{eq:power-euclidean} by Theorem~\ref{thm:hms}. In one dimension
	the sign is constant on each half-line, giving the same derivative
	estimate for $\sgn(t)|t|^{is}$.
	
	To justify the substitutions, including on sets where the arguments vanish,
	we use the continuous regularizations
	\begin{equation}\label{eq:regularization}
		\Psi_{s,\varepsilon}(t):=t(t^2+\varepsilon^2)^{(is-1)/2},
		\qquad
		\Phi_{n,s,\varepsilon}(z)
		:=|z|^2(|z|^2+\varepsilon^2)^{is/2-1}.
	\end{equation}
	They vanish at zero and converge pointwise everywhere to $\Psi_s$ and
	$\Phi_{n,s}$, respectively. Their Mikhlin quantities satisfy, uniformly
	in $\varepsilon>0$,
	\[
	\mathfrak M_1(\Psi_{s,\varepsilon})\lesssim1+|s|,\qquad
	\mathfrak M_{\lfloor n/2\rfloor+1}(\Phi_{n,s,\varepsilon})
	\lesssim_n(1+|s|)^{\lfloor n/2\rfloor+1}.
	\]
	Scaling $z=\varepsilon w$ reduces the second assertion
	to the function $|w|^2(1+|w|^2)^{is/2-1}$, up to the unimodular
	factor $\varepsilon^{is}$. The Leibniz rule shows that each derivative
	of order $k$, multiplied by $|w|^k$, is bounded by
	$C_k(1+|s|)^k$: for $|w|\leq1$ all remaining powers are bounded,
	and for $|w|\geq1$ the powers of $|w|$ cancel. The same argument
	applies to $w(1+w^2)^{(is-1)/2}$.
	
	We apply Theorem~\ref{thm:hms} to these regularizations.
	Lemma~\ref{lem:pullback}, with $\alpha=\beta=b$, proves
	\eqref{eq:power-signed}. With $\alpha=g$ and $\beta=-g$ it proves
	\eqref{eq:power-sum}, since $\alpha(\xi)-\beta(\eta)=g(\xi)+g(\eta)$.
\end{proof}

In the next section, we will need the following result about the Riesz difference multiplier, defined for $1\leq j\leq d$ by
\[
R_j(\xi,\eta):=
\begin{cases}(\xi_j-\eta_j)/|\xi-\eta|,&\xi\ne\eta,\\
	0,&\xi=\eta.\end{cases}
\]

\begin{lemma}\label{lem:riesz}
	For $1<p<\infty$ and $1\leq j\leq d$, we have
	\[
	\|\Sch_{R_j}\|_{\Sp^p\to\Sp^p}\lesssim_{d,p}1 \,,
	\qquad
	\|\Sch_{R_j}\|_{\Sp^p_{\rm weak}\to\Sp^p_{\rm weak}}\lesssim_{d,p}1 \,.
	\]
\end{lemma}

\begin{proof}
	The function $\rho_j(z):=z_j/|z|$ is smooth away from zero and homogeneous
	of degree zero. Hence $\partial^\alpha\rho_j$ is homogeneous of degree
	$-|\alpha|$, and its restriction to the unit sphere is bounded.
	Thus $\mathfrak M_{\lfloor d/2\rfloor+1}(\rho_j)\lesssim_d1$ and
	Theorem~\ref{thm:hms} proves the strong estimate. Real interpolation
	between $\Sp^{p_0}$ and $\Sp^{p_1}$ with $1<p_0<p<p_1<\infty$
	gives the weak estimate.
\end{proof}

\begin{proof}[Proof of Proposition~\ref{prop:weighted-multiplier}]
It suffices to consider real-valued $b$. Indeed, both $\operatorname{Re}b$
and $\operatorname{Im}b$ satisfy the same hypothesis with the same $g$,
and the denominator of $r_\theta$ is independent of $b$.
The assertion for complex $b$ follows by linearity and the triangle
or quasitriangle inequality in the corresponding Schatten class.

By symmetry, we may assume that $1/2\leq\theta\leq1$. We put
\[
D(\xi,\eta):=g(\xi)^\theta g(\eta)^{1-\theta}
             +g(\xi)^{1-\theta}g(\eta)^\theta.
\]
If $\theta<1$ and $\{g=0\}$ has positive measure, then
$D(\xi,\eta)=0$ whenever $g(\xi)=0$. By
\eqref{eq:weighted-hajlasz} and Fubini's theorem, $b$ is constant
almost everywhere, and $r_\theta=0$. We may therefore assume that
$g>0$ almost everywhere when $\theta<1$.
For $s\in\R$, we have
\[
\Phi_{1,s}(D(\xi,\eta))
=g(\xi)^{is(1-\theta)}g(\eta)^{is(1-\theta)}
 \Phi_{1,s}\bigl(g(\xi)^{2\theta-1}+g(\eta)^{2\theta-1}\bigr),
\]
where the first two factors are interpreted as one when $\theta=1$,
including where $g=0$. They act by left and right multiplication by
unitaries, while \eqref{eq:power-sum}, applied to $g^{2\theta-1}$,
bounds the last factor. Thus, uniformly in $\theta$,
\[
\|\Sch_{\Phi_{1,s}(D)}\|_{\Sp^p\to\Sp^p}\lesssim_p1+|s|.
\]

We choose $\chi\in C^\infty(\R)$ equal to one on $(-\infty,0]$ and
zero on $[1,\infty)$, and set $H(t):=e^t\chi(t)$. Then
$H\in\mathcal S(\R)$, and with
\[
c(s):=\frac1{2\pi}\int_\R H(t)e^{-ist}\,dt,
\qquad H(t)=\int_\R c(s)e^{ist}\,ds,
\]
we have $c\in\mathcal S(\R)$. We put $\delta:=b(\xi)-b(\eta)$ and
$A:=|\xi-\eta|D(\xi,\eta)$. If $\delta\ne0$, then
$0<|\delta|/A\leq1$ and
$\delta/A=\sgn(\delta)H(\log(|\delta|/A))$. Consequently,
\begin{equation}\label{eq:mellin}
r_\theta(\xi,\eta)=\int_\R c(s)
\Psi_s(b(\xi)-b(\eta))\Phi_{d,-s}(\xi-\eta)
\Phi_{1,-s}(D(\xi,\eta))\,ds.
\end{equation}
If $\delta=0$, both sides vanish; a zero denominator also forces
$\delta=0$ by \eqref{eq:weighted-hajlasz}.
Products of Schur kernels correspond to compositions of their multipliers.
Lemma~\ref{lem:powers} and the bound above give
\[
\|\Sch_{r_\theta}\|_{\Sp^p\to\Sp^p}
\lesssim_{d,p}\int_\R |c(s)|(1+|s|)^{\lfloor d/2\rfloor+3}\,ds<\infty.
\]
More precisely, on finite-rank inputs the integrand is continuous in
$\Sp^2$ by dominated convergence. Its pairing with a finite-rank
$\Sp^{p'}$ operator is bounded by the displayed integrable majorant.
Integrating these pairings and using Schatten duality proves the
$\Sp^p$ estimate, first on finite-rank inputs and then by density.
Finally, real interpolation between two exponents
$1<p_0<p<p_1<\infty$ proves the weak estimate.
\end{proof}

%%%%%%%%%%%%%%%%%%%%%%%%

\section{The commutator identity and proof of Schatten class estimates}\label{sec:commid}

In this section we prove Theorems \ref{thm:main} and \ref{thm:bscomm}. In Subsection \ref{sec:definitions} we define product operators and commutators as distributional forms. Then, in Subsection \ref{sec:commidsubsec} we prove a crucial commutator identity involving sandwiched products. After some technical preliminaries in Subsection \ref{sec:hajlasz}, in Subsection \ref{sec:proofmain} we combine the commutator identity with classical estimates for $f(X)g(-\ii\nabla)$ to derive our main results.

%%%%%%%%%%%%%%%%%%%%%%%%

\subsection{Product and commutator forms}\label{sec:definitions}

We use the unitary Fourier transform
\[
\F\psi(\xi)=\widehat\psi(\xi)
:=(2\pi)^{-d/2}\int_{\R^d}e^{-\ii x\cdot\xi}\psi(x)\dd x,
\]
and take inner products to be linear in the second entry.
Our test functions belong to the Schwartz space $\mathcal S(\R^d)$,
which is preserved by the Fourier transform.
We first recall how to define the product operators $f(X)g(-\ii\nabla)$ used below.
We suppose that $f,g\in L^2_{\loc}(\R^d)$ and that, for some $N\geq0$,
\[
\sup_{k\in\Z^d}(1+|k|)^{-N}
\bigl(\|f\|_{L^2(Q_k)}+\|g\|_{L^2(Q_k)}\bigr)<\infty.
\]
Here, as everywhere, $Q_k=k+[0,1)^d$ for $k\in\Z^d$. For $\phi,\psi\in\mathcal S(\R^d)$, we set
\begin{equation}\label{eq:product-forms}
\begin{split}
\langle\phi,f(X)g(-\ii\nabla)\psi\rangle
&:=\langle\overline f\,\phi,g(-\ii\nabla)\psi\rangle,\\
\langle\phi,g(-\ii\nabla)f(X)\psi\rangle
&:=\langle\overline g(-\ii\nabla)\phi,f\psi\rangle.
\end{split}
\end{equation}
These are ordinary $L^2$ pairings. Indeed, the rapid decay of the test
functions and their Fourier transforms gives
\begin{equation}\label{eq:gpsi}
\begin{split}
\|f\psi\|_2^2
&\leq\sum_k\|f\|_{L^2(Q_k)}^2
            \sup_{x\in Q_k}|\psi(x)|^2<\infty,\\
\|g\widehat\psi\|_2^2
&\leq\sum_k\|g\|_{L^2(Q_k)}^2
            \sup_{\xi\in Q_k}|\widehat\psi(\xi)|^2<\infty.
\end{split}
\end{equation}
Thus $f\psi\in L^2$ and $g(-\ii\nabla)\psi=\F^*(g\widehat\psi)\in L^2$.
The same argument shows that $\overline f\,\phi,\overline g(-\ii\nabla)\phi\in L^2$.
The resulting forms are continuous on pairs of Schwartz functions and agree
with the operator products whenever those are defined. Saying that such
a form belongs to a Schatten class means that it is represented by a
bounded operator in that class; this operator is unique by density.

We next verify the corresponding definition for the commutator under
a common hypothesis satisfied by all the gradient spaces considered here.

\begin{lemma}\label{prop:comm-form}
	We assume that $a,b\in L^1_{\rm loc}(\R^d)$ have gradients satisfying, for some
	$\max\{1,2d/(d+2)\}<r<\infty$,
	\[
	\sup_k\|\nabla a\|_{L^r(Q_k)}<\infty,
	\qquad \sup_k\|\nabla b\|_{L^r(Q_k)}<\infty.
	\]
	For $\phi,\psi\in\mathcal S(\R^d)$ the formula
	\begin{equation}\label{eq:comm-form}
		\langle\phi,[a(X),b(-\ii\nabla)]\psi\rangle
		:=\langle\overline a\,\phi,b(-\ii\nabla)\psi\rangle
		-\langle\overline b(-\ii\nabla)\phi,a\psi\rangle
	\end{equation}
	is well defined. It is a continuous sesquilinear form on pairs of
	Schwartz functions, agrees with the usual commutator wherever the operator
	products are defined, and is unchanged by adding constants to either
	function.
\end{lemma}

\begin{proof}
	We apply Lemma~\ref{lem:unit-growth} to both $a$ and $b$. It gives
	$a,b\in L^2_{\rm loc}$ and
	\[
	\|a\|_{L^2(Q_k)}+\|b\|_{L^2(Q_k)}\lesssim_{a,b}1+|k|.
	\]
	Consequently both product forms in \eqref{eq:product-forms}, with
	$f=a$ and $g=b$, are well defined and continuous. Their difference is
	\eqref{eq:comm-form}.
	Moving the left factor of each operator product to the first argument
	of the inner product shows agreement with the usual commutator whenever
	the products are defined. Cancellation of a constant added to either function follows easily from this definition.
\end{proof}

\begin{remark}\label{rem:comm-form}
The gradient assumptions of Lemma~\ref{prop:comm-form} are not needed
for the definition itself. Indeed, the discussion of product forms
shows that \eqref{eq:comm-form} is well defined and continuous on pairs
of Schwartz functions $\phi,\psi\in\mathcal S(\R^d)$ whenever
$a,b\in L^2_{\loc}$ and the $L^2$ norms of $a$ and $b$ on the cubes $Q_k$ grow
at most polynomially.
\end{remark}

%%%%%%%%%%%%%%%%%%%%%

\subsection{The commutator identity}\label{sec:commidsubsec}

The following proposition is the key ingredient in the proofs of Theorems \ref{thm:main} and \ref{thm:bscomm}. It expresses the commutator in terms of Schur multipliers applied to sandwiched products. We state the identity for smooth $a$ and pass to general $a$ in the proofs of the main results.

\begin{proposition}\label{prop:commutator-identity}
For $1<p<\infty$ and $\mathcal J$ equal to either $\Sp^p$ or
$\Sp^p_{\rm weak}$, we assume that $a\in C_c^\infty(\R^d)$,
$b\in L^2_{\loc}(\R^d)$, and that, for some $N\geq 0$,
\[
\sup_k(1+|k|)^{-N}\|b\|_{L^2(Q_k)}<\infty.
\]
We further assume that $0\leq\theta\leq1$, and
that a measurable $g:\R^d\to[0,\infty)$ satisfies
\eqref{eq:weighted-hajlasz} and
\[
\sup_k\bigl(\|g^\theta\|_{L^2(Q_k)}
            +\|g^{1-\theta}\|_{L^2(Q_k)}\bigr)<\infty.
\]
We assume that, for each $j=1,\ldots,d$, the form
\[
K_j:=g^\theta(-\ii\nabla)(\partial_j a)(X)g^{1-\theta}(-\ii\nabla)
   +g^{1-\theta}(-\ii\nabla)(\partial_j a)(X)g^\theta(-\ii\nabla)
\]
extends to an operator in $\mathcal J$. We set
\[
m_j(\xi,\eta):=
\frac{(b(\xi)-b(\eta))(\xi_j-\eta_j)}
{|\xi-\eta|^2\bigl(g(\xi)^\theta g(\eta)^{1-\theta}
                  +g(\xi)^{1-\theta} g(\eta)^\theta\bigr)},
\]
with value zero when the denominator vanishes. Then the form
$[a(X),b(-\ii\nabla)]$ extends to an operator in $\mathcal J$, and
\begin{equation}\label{eq:weighted-identity}
\F[a(X),b(-\ii\nabla)]\F^*
=\ii\sum_{j=1}^d\Sch_{m_j}(\F K_j\F^*)
\end{equation}
and
\[
\|[a(X),b(-\ii\nabla)]\|_{\mathcal J}
\lesssim_{d,p}\sum_{j=1}^d\|K_j\|_{\mathcal J}.
\]
\end{proposition}

We define $K_j$ as a sesquilinear form by moving the left Fourier
multiplier in each term to the first argument of the inner product.
This form is well defined and continuous on pairs of Schwartz
functions, since $\partial_j a$ is bounded and \eqref{eq:gpsi} applies
with $g$ replaced by $g^\theta$ or $g^{1-\theta}$.
The commutator form is well defined and continuous on the same class
by Remark~\ref{rem:comm-form}.

\begin{proof}
Proposition~\ref{prop:weighted-multiplier} and
Lemma~\ref{lem:riesz} give
\[
\|\Sch_{m_j}\|_{\mathcal J\to\mathcal J}\lesssim_{d,p}1.
\]

For $R>0$, let $P_R:=\1_{B_R(0)}(-\ii\nabla)$. Thus
$\F P_R\F^*$ is multiplication by $\1_{B_R(0)}$.
The operator $\F P_R K_jP_R\F^*$ has integral kernel zero
outside $B_R(0)\times B_R(0)$ and equal to
\[
(2\pi)^{-d/2}\ii(\xi_j-\eta_j)\widehat a(\xi-\eta)
\bigl(g(\xi)^\theta g(\eta)^{1-\theta}
     +g(\xi)^{1-\theta} g(\eta)^\theta\bigr)
\]
on that set. Indeed, the formula follows by Fubini's theorem from the
definition of the forms, first for
$\widehat\phi,\widehat\psi\in C_c^\infty(B_R(0))$.
Since $\widehat{\partial_j a}$ is bounded and
$g^\theta,g^{1-\theta}\in L^2(B_R(0))$, this kernel belongs to $L^2$.
The defining action of Schur multiplication gives
\[
 \Sch_{m_j}(\F P_R K_jP_R\F^*)
 =(\F P_R\F^*)\Sch_{m_j}(\F K_j\F^*)(\F P_R\F^*).
\]
Multiplying the displayed kernel by $m_j$, summing over $j$, and
multiplying by $\ii$ gives, on $B_R(0)\times B_R(0)$,
\[
(2\pi)^{-d/2}\bigl(b(\eta)-b(\xi)\bigr)\widehat a(\xi-\eta),
\]
where we also used \eqref{eq:weighted-hajlasz} to conclude that $b(\xi)=b(\eta)$ almost everywhere on the set where the denominator vanishes. This kernel also represents the commutator form when both Fourier transforms
are supported in $B_R(0)$. Thus \eqref{eq:weighted-identity} holds for
Schwartz functions with smooth compactly supported Fourier transforms.
These functions are dense in the Schwartz class, so continuity proves
the identity on all pairs of Schwartz functions.
Its right side belongs to $\mathcal J$ and
therefore gives the asserted extension and estimate.
\end{proof}

%%%%%%%%%%%%%%%%%%%%%%%%%%%%%

\subsection{A Haj\l asz gradient bound}\label{sec:hajlasz}

The following lemma supplies the functions $g$ needed in the inequality \eqref{eq:weighted-hajlasz}. It involves the Hardy--Littlewood maximal operator $M$, which is defined for  $f\in L^1_{\rm loc}(\R^d)$ by
\[
Mf(x):=\sup_{r>0}\frac{1}{|B_r(x)|}\int_{B_r(x)}|f(y)|\,dy.
\]

\begin{lemma}\label{lem:hajlasz-gradients}
For $b\in W^{1,1}_{\loc}(\R^d)$, the following assertions hold.
\begin{enumerate}
\item[(i)] With $G:=M|\nabla b|$, one has
\[
|b(\xi)-b(\eta)|\lesssim_d|\xi-\eta|
\bigl(G(\xi)+G(\eta)\bigr).
\]
\item[(ii)] If in addition $d\geq2$, $1<q<d$, and
$b\in W^{1,q}_{\loc}(\R^d)$, then, with
$G:=(M(|\nabla b|^q))^{1/q}$, one has
\begin{equation}\label{eq:refined-hajlasz}
\begin{split}
|b(\xi)-b(\eta)|\lesssim_{d,q}|\xi-\eta|
\bigl(&G(\xi)^{1-q/d}G(\eta)^{q/d}+G(\xi)^{q/d}G(\eta)^{1-q/d}\bigr).
\end{split}
\end{equation}
\item[(iii)] If $d=1$, then, with $G:=M|b'|$,
\begin{equation}\label{eq:refined-hajlasz-one}
|b(\xi)-b(\eta)|\leq2|\xi-\eta|
\min\{G(\xi),G(\eta)\}
\leq2|\xi-\eta|G(\xi)^{1/2}G(\eta)^{1/2}.
\end{equation}
\end{enumerate}
Each inequality holds at all Lebesgue points $\xi,\eta$ of $b$.
\end{lemma}

Part (i) is well known; see \cite[Section~2]{Ha}. Its proof
is included for completeness. The refinement in part (ii) follows by
optimizing the estimates for ball averages, while part (iii) is an
elementary one-dimensional observation.

\begin{proof}
	For (i), the inequality is immediate if either value of $G$ is infinite. We therefore consider Lebesgue points $\xi,\eta$ of $b$ at which $G$ is finite, with \(r:=|\xi-\eta|>0\) and set $\zeta:=(\xi+\eta)/2$. Since \(\xi,\eta\in B_r(\zeta)\), the version of the bound \eqref{eq:poincare-potential} for balls gives
	\[
	\bigl|b(\xi)-b_{B_r(\zeta)}\bigr|
	\lesssim_d
	\int_{B_r(\zeta)}
	\frac{|\nabla b(\omega)|}{|\xi-\omega|^{d-1}}\,d\omega.
	\]
	Now \(B_r(\zeta)\subset B_{3r/2}(\xi)\), and the elementary estimate
	\[
	\int_{B_R(\xi)}
	\frac{|f(\omega)|}{|\xi-\omega|^{d-1}}\,d\omega
	\le d|B_1(0)|\,R\,Mf(\xi) \,,
	\]
	shows that
	\[
	\bigl|b(\xi)-b_{B_r(\zeta)}\bigr|
	\lesssim_d r\,M(|\nabla b|)(\xi).
	\]
	The same argument for $\eta$ instead of $\xi$ yields
	\[
	\begin{aligned}
		|b(\xi)-b(\eta)|
		&\le
		\bigl|b(\xi)-b_{B_r(\zeta)}\bigr|
		+\bigl|b(\eta)-b_{B_r(\zeta)}\bigr|
		\lesssim_d r\bigl(M(|\nabla b|)(\xi)+M(|\nabla b|)(\eta)\bigr) \,.
	\end{aligned}
	\]
	This proves (i).

For (ii) and (iii), if the gradient vanishes almost everywhere, $b$ is constant and the
assertions are immediate. Otherwise $G$ is strictly positive everywhere.
For (ii), the inequality is immediate if either value of $G$ is infinite,
so we assume both values are finite. We put $r:=|\xi-\eta|>0$ and suppose that
$G(\xi)\geq G(\eta)$. For $0<t\leq r$, H\"older's inequality and
$B_t(\xi)\subset B_{2r}(\eta)$ give
\[
\frac{1}{|B_t(\xi)|}\int_{B_t(\xi)}|\nabla b(\omega)|\,d\omega
\lesssim_d\min\left\{G(\xi),(r/t)^{d/q}G(\eta)\right\}.
\]
Telescoping the averages on the balls $B_{2^{-k}r}(\xi)$, using the
$L^1$ Poincar\'e inequality and the Lebesgue-point property, yields
\[
\begin{split}
|b(\xi)-b_{B_r(\xi)}|
&\lesssim_d\sum_{k\geq0}2^{-k}r
\min\{G(\xi),2^{kd/q}G(\eta)\}\\
&\lesssim_{d,q}r\,G(\xi)^{1-q/d}G(\eta)^{q/d}.
\end{split}
\]
For the last bound, we split the sum where $2^{-k}$ is comparable to
$(G(\eta)/G(\xi))^{q/d}$. Below that radius we use $G(\xi)$; above it
we use $2^{kd/q}G(\eta)$. Both geometric sums have the displayed bound,
since $q<d$.

The same telescoping argument at $\eta$ gives
$|b(\eta)-b_{B_r(\eta)}|\lesssim_d rG(\eta)$.
Poincar\'e on $B_{3r}(\eta)$, which contains both balls of radius $r$,
also gives
\[
|b_{B_r(\xi)}-b_{B_r(\eta)}|\lesssim_d rG(\eta).
\]
Since $G(\eta)\leq G(\xi)^{1-q/d}G(\eta)^{q/d}$, these estimates
prove \eqref{eq:refined-hajlasz}. The opposite ordering of the two
values of $G$ is treated by interchanging the points. Thus the only excluded points
are the non-Lebesgue points, independently of the balls used in the proof.

For (iii), we put $r:=|\xi-\eta|>0$ and integrate $|b'|$ over the interval between the two
points. This interval lies in each of $B_r(\xi)$ and $B_r(\eta)$, so
its integral is at most $2r\min\{G(\xi),G(\eta)\}$. This proves
\eqref{eq:refined-hajlasz-one}.
\end{proof}

%%%%%%%%%%%%%%%%%%%%%%%%%%%%%

\subsection{Proof of the commutator bounds}\label{sec:proofmain}

With the commutator identity in Proposition \ref{prop:commutator-identity}, the Schur multiplier bound in Proposition \ref{prop:weighted-multiplier} and the gradient bounds in Lemma \ref{lem:hajlasz-gradients} at our disposal, it is straightforward to prove our positive results in Theorems \ref{thm:main} and \ref{thm:bscomm}.

We begin with the Lebesgue-space bounds in Theorem~\ref{thm:main}. In \eqref{eq:weighted-identity}, the choice $\theta=1$ suffices for part (a) when $p>2$ and for part (b) when $p\geq2$. To reach $p\leq2$ in part (a), or $p<2$ in part (b), we use $0<\theta<1$.

\begin{proof}[Proof of Theorem~\ref{thm:main}]
We set $E:=L^p_{\rm weak}$ in part (a) and $E:=L^p$ in part (b).
We assume that $b\in L^1_{\loc}(\R^d)$ satisfies $\nabla b\in E$.

In dimension $d\geq2$ we choose
\begin{equation}\label{eq:choice-q}
\max\{1,d(1-p/2)\}<q<\min\{p,d\},
\qquad \theta:=1-q/d.
\end{equation}
Such a choice is possible because $p>1$ and $p>2d/(d+2)$. It gives
\begin{equation}\label{eq:factor-exponents}
0<\theta<1,\qquad
p/\theta>2,\quad p/(1-\theta)>2.
\end{equation}
In particular, $p/\theta>2$ follows from $q>d(1-p/2)$, while
$p/(1-\theta)=dp/q>d\geq2$.
Since $q<p$, the gradient assumption gives $b\in W^{1,q}_{\loc}$. By Lemma~\ref{lem:hajlasz-gradients}(ii), a sufficiently large constant
$C_{d,q}$ makes \eqref{eq:weighted-hajlasz} hold with
\begin{equation}\label{eq:g-refined}
g:=C_{d,q}\bigl(M(|\nabla b|^q)\bigr)^{1/q}.
\end{equation}
In dimension one we take $\theta:=1/2$ and $g:=M|b'|$ instead, using Lemma~\ref{lem:hajlasz-gradients} (iii).
Then \eqref{eq:weighted-hajlasz} and \eqref{eq:factor-exponents} again hold.

The strong and weak maximal theorems give
\begin{equation}\label{eq:g-refined-norm}
\|g\|_E\lesssim_{d,p}\|\nabla b\|_E.
\end{equation}
By \eqref{eq:factor-exponents}, both $g^\theta$ and $g^{1-\theta}$ have
uniformly bounded $L^2$ norms on unit cubes.
Lemma~\ref{lem:unit-growth} also gives $b\in L^2_{\loc}$ and
$\|b\|_{L^2(Q_k)}\lesssim_b1+|k|$: we choose a local gradient exponent
strictly between $\max\{1,2d/(d+2)\}$ and $p$.

We first assume that $a\in C_c^\infty(\R^d)$. For $f:=\partial_j a$,
we factor
\begin{equation}\label{eq:sandwich-factorization}
\begin{split}
&g^\theta(-\ii\nabla)f(X)g^{1-\theta}(-\ii\nabla)=\bigl(g^\theta(-\ii\nabla)|f|^\theta(X)\bigr)
\sgn(f)(X)
\bigl(|f|^{1-\theta}(X)g^{1-\theta}(-\ii\nabla)\bigr),
\end{split}
\end{equation}
where $\sgn(f):=f/|f|$ when $f\ne0$, and zero otherwise. (A similar asymmetric factorization also appears in \cite[Lemma~4.7]{FSZ}.) For the weak estimate we apply Cwikel's inequality at exponents
$p/\theta$ and $p/(1-\theta)$, and for the strong estimate we apply the
Kato--Seiler--Simon inequality at the same exponents; see
\cite[Theorems~4.1 and~4.2]{Simon}.
The two factor norms are bounded by
\[
\bigl(\|f\|_{L^p}\|g\|_E\bigr)^\theta
\quad\text{and}\quad
\bigl(\|f\|_{L^p}\|g\|_E\bigr)^{1-\theta},
\]
up to constants depending only on $d,p$. H\"older's inequality for
strong or weak Schatten classes, with
$\theta/p+(1-\theta)/p=1/p$, therefore gives
\begin{equation}\label{eq:sandwich-bound}
\|g^\theta(-\ii\nabla)f(X)g^{1-\theta}(-\ii\nabla)\|_{\mathcal J}
\lesssim_{d,p}\|f\|_{L^p}\|g\|_E,
\end{equation}
where $\mathcal J$ is the corresponding strong or weak Schatten class.
The factorization defines a bounded operator and agrees with the form
in Proposition~\ref{prop:commutator-identity}: we move the first factor to the left
in the $L^2$ pairing and use
$|f|^\theta\sgn(f)|f|^{1-\theta}=f$.
The same estimate holds with $\theta$ replaced by $1-\theta$.
Proposition~\ref{prop:commutator-identity} and \eqref{eq:g-refined-norm} now prove
both parts of the theorem for smooth $a$.

For general $a\in L^1_{\rm loc}(\R^d)$ with $\nabla a\in L^p(\R^d)$, Lemma~\ref{lem:approximation} supplies
$a_n\in C_c^\infty$ and constants $c_n$ such that
\[
\nabla a_n\longrightarrow\nabla a\quad\text{in }L^p,
\qquad a_n+c_n\longrightarrow a\quad\text{in }L^2_{\loc}.
\]
Applying the smooth estimate to $a_n-a_m$ shows that the commutators
are Cauchy in $\mathcal J$, which is complete. To identify the limit, it suffices
to take $\phi,\psi\in C_c^\infty$, since $C_c^\infty$ is dense in
$\mathcal S$ and both the commutator form and the form of the limiting
bounded operator are continuous on pairs of Schwartz functions.
Constants cancel in the commutator form, and
\[
\begin{split}
&|\langle\phi,([a_n(X),b(-\ii\nabla)]-[a(X),b(-\ii\nabla)])\psi\rangle|\\
&\quad\leq\|(a_n+c_n-a)\phi\|_2\|b(-\ii\nabla)\psi\|_2 +\|\overline b(-\ii\nabla)\phi\|_2\|(a_n+c_n-a)\psi\|_2
\longrightarrow0.
\end{split}
\]
Thus the limit represents the commutator form \eqref{eq:comm-form}.
Passing to the limit proves the estimates.
\end{proof}

Now we prove the bounds in Birman--Solomyak spaces in Theorem \ref{thm:bscomm} and emphasize that in this case the commutator identity \eqref{eq:weighted-identity} with $\theta=1$ suffices.

\begin{proof}[Proof of Theorem~\ref{thm:bscomm}]
We set $E:=\ell^p_{\rm weak}(L^2)$ in part (a) and
$E:=\ell^p(L^2)$ in part (b).
We assume that $b\in L^1_{\loc}(\R^d)$ satisfies $\nabla b\in E$.

By Lemma~\ref{lem:hajlasz-gradients}(i), \eqref{eq:hajlasz} holds with
\begin{equation}\label{eq:g}
g:=C_dM(|\nabla b|).
\end{equation}
The maximal function bound in Lemma~\ref{lem:bsmax} below gives
\begin{equation}\label{eq:g-norm}
\|g\|_E\lesssim_{d,p}\|\nabla b\|_E.
\end{equation}
In particular, $g$ has uniformly bounded $L^2$ norms on unit cubes.
Lemma~\ref{lem:unit-growth} gives $b\in L^2_{\loc}$ and
$\|b\|_{L^2(Q_k)}\lesssim_b1+|k|$.

We first assume that $a\in C_c^\infty(\R^d)$. We put $f:=|\nabla a|$ and
$T:=f(X)g(-\ii\nabla)$. The Birman--Solomyak--Simon estimates
\cite[Theorems~4.5 and~4.6]{Simon} yield
\[
\|T\|_{\mathcal J}
\lesssim_{d,p}\|\nabla a\|_{\ell^p(L^2)}\|g\|_E,
\]
where $\mathcal J$ is the corresponding strong or weak class.
We put $h_j:=(\partial_j a)/f$, with value zero where $f=0$.
Since $|h_j|\leq1$ and $f,g$ are real, the products in
Proposition~\ref{prop:commutator-identity}, with
$\theta=1$, satisfy
\[
K_j=T^*h_j(X)+h_j(X)T,
\qquad \|K_j\|_{\mathcal J}\lesssim_p\|T\|_{\mathcal J}.
\]
These equalities hold as forms and provide the required bounded
extensions. The proposition and \eqref{eq:g-norm} prove both estimates
for smooth $a$.

For general $a\in L^1_{\rm loc}$ with $\nabla a\in\ell^p(L^2)$, Lemma~\ref{lem:approximation} supplies
$a_n\in C_c^\infty$ and constants $c_n$ such that
\[
\nabla a_n\longrightarrow\nabla a\quad\text{in }\ell^p(L^2),
\qquad a_n+c_n\longrightarrow a\quad\text{in }L^2_{\loc}.
\]
We keep $b$ fixed and apply the smooth estimate to $a_n-a_m$.
The commutators are Cauchy in $\mathcal J$, which is complete. Local $L^2$ convergence and cancellation
of constants identify their limit with the commutator form, as in the
proof of Theorem~\ref{thm:main}.
Passing to the limit proves the assertions.
\end{proof}

\begin{remark}\label{rem:simon}
	At $p=2$ there is a direct argument, essentially due to \cite{Simon79}, which gives the constant in \eqref{eq:strong2}. We first assume that $a,b\in C_c^\infty(\R^d)$. We have
	\begin{align*}
		\|[a(X),b(-\ii\nabla)]\|_{\Sp^2}^2
		&=(2\pi)^{-d}\int_{\R^d}|\widehat a(z)|^2
		\int_{\R^d}|b(\eta+z)-b(\eta)|^2\,d\eta\,dz\\
		&\leq (2\pi)^{-d}\int_{\R^d}|z|^2|\widehat a(z)|^2\,dz
		\|\nabla b\|_2^2\\
		&=(2\pi)^{-d}\|\nabla a\|_2^2\|\nabla b\|_2^2.
	\end{align*}
	The middle inequality follows from
	\[
	b(\eta+z)-b(\eta)
	=\int_0^1 z\cdot\nabla b(\eta+tz)\,dt
	\]
	and Minkowski's inequality. For general $a,b$, Lemma~\ref{lem:approximation}
	gives $a_n,b_n\in C_c^\infty$ with
	$\nabla a_n\to\nabla a$ and $\nabla b_n\to\nabla b$ in $L^2$.
	The already proved estimate in Theorem~\ref{thm:main}(b), with $p=2$,
	gives
	\[
	\begin{split}
	&\|[a_n(X),b_n(-\ii\nabla)]-[a(X),b(-\ii\nabla)]\|_{\Sp^2}\\
	&\quad\lesssim_d
	\|\nabla(a_n-a)\|_2\|\nabla b_n\|_2
	+\|\nabla a\|_2\|\nabla(b_n-b)\|_2
	\longrightarrow0.
	\end{split}
	\]
	Passing to the limit in the smooth estimate therefore gives the
	same explicit constant for $a,b$.
\end{remark}

The proof of Theorem~\ref{thm:bscomm} used the boundedness of the Hardy--Littlewood maximal operator on $\ell^p(L^2)$ and $\ell^p_{\rm weak}(L^2)$. For the sake of completeness, let us provide a proof.

\begin{lemma}\label{lem:bsmax}
	For $1<p<\infty$ and $f\in\ell^p(L^2)$, respectively $f\in\ell^p_{\rm weak}(L^2)$, we have
	\[
	\|Mf\|_{\ell^p(L^2)}\lesssim_{d,p}\|f\|_{\ell^p(L^2)},\qquad
	\|Mf\|_{\ell^p_{\rm weak}(L^2)}
	\lesssim_{d,p}\|f\|_{\ell^p_{\rm weak}(L^2)}.
	\]
\end{lemma}

\begin{proof}
We put $f_j:=\|f\|_{L^2(Q_j)}$ and
$F:=\sum_j f_j\mathbf1_{Q_j}$. We split the balls centered at a point
of $Q_k$ into those of radius at most $2\sqrt d$ and the remaining
balls. The first contribution is bounded by the maximal function of
$f$ restricted to a fixed finite union of neighboring cubes. The
$L^2$ maximal theorem bounds its norm on $Q_k$ by
$C_d(\sum_{|j-k|\leq C_d}f_j^2)^{1/2}$.

For $x,y\in Q_k$ and $r>2\sqrt d$, every unit cube meeting $B_r(x)$
is contained in $B_{2r}(y)$. Cauchy--Schwarz on each unit cube gives
\[
\frac1{|B_r(x)|}\int_{B_r(x)}|f(z)|\,dz
\leq\frac1{|B_r(x)|}\sum_{Q_j\cap B_r(x)\ne\varnothing}f_j
\leq\frac1{|B_r(x)|}\int_{B_{2r}(y)}F(z)\,dz
\lesssim_d MF(y).
\]
Taking the infimum over $y\in Q_k$, we obtain
\[
\|Mf\|_{L^2(Q_k)}
\lesssim_d\left(\sum_{|j-k|\leq C_d}f_j^2\right)^{1/2}
+\inf_{y\in Q_k}MF(y).
\]
Finite sums of translates are bounded on both sequence spaces in
the statement. Moreover, the piecewise constant function with value
$\inf_{Q_k}MF$ on $Q_k$ is bounded by $MF$. The strong and weak
maximal theorems on $\R^d$ therefore bound the corresponding sequence
norms, up to constants depending only on $d,p$, by
$\|F\|_{L^p}=\|(f_j)\|_{\ell^p}$ and
$\|F\|_{L^p_{\rm weak}}=\|(f_j)\|_{\ell^p_{\rm weak}}$, respectively.
This proves both assertions.
\end{proof}

%%%%%%%%%%%%%%%%%%%%%%%%%%%%%%%%
%%%%%%%%%%%%%%%%%%%%%%%%%%%%%%%%%%%%%%

\section{Failure of the Lipschitz endpoint estimate}\label{sec:endpoint}

In this section we prove Theorem~\ref{thm:counterexample}. We first construct real smooth
bounded periodic functions on $\R$ with derivative norms at most one and
arbitrarily large commutator norms. Tensor products and smooth cutoffs
then give compactly supported functions in every dimension whose
commutator norms are arbitrarily large while both gradient norms
are at most one.

The construction uses position frequencies $R^j$ and momentum frequencies
$2\pi R^{-k}$ for a large integer $R$. We construct oscillating functions with these frequencies and take their primitives. In the commutator of these primitives, the interactions with $j<k$ survive as $R\to\infty$, while the others vanish. Lemma \ref{lem:polynomials} constructs bounded trigonometric polynomials whose surviving
interactions have a large sum; Lemma \ref{lem:norm} turns that sum into an
operator-norm lower bound.

%%%%%%%%%%%%%%%%%%%%%%%%%%%%

\subsection{Two bounded trigonometric polynomials}

A trigonometric polynomial on $\T^N$, where $\T:=\R/(2\pi\Z)$, is a finite sum
\[
H(z)=\sum_{m\in\Z^N}\widehat H(m)e^{\ii m\cdot z}.
\]
We call $e^{\ii m\cdot z}$ a \emph{Fourier monomial} and the integers
$m_1,\ldots,m_N$ its \emph{coordinate Fourier exponents}. A monomial
occurs in $H$ if its coefficient $\widehat H(m)$ is nonzero.

\begin{lemma}\label{lem:polynomials}
	For every integer $N\ge2$, there are real trigonometric
	polynomials $F,G$ on $\T^N$, with decompositions into real trigonometric
	polynomials
	\[
	F=\sum_{j=1}^N F_j,\qquad G=\sum_{k=1}^N G_k,
	\]
	having the following properties:
	\begin{enumerate}
		\item[(i)] $\|F\|_\infty\le1$ and $\|G\|_\infty\le1$.
		\item[(ii)] All coordinate Fourier exponents belong to $\{-1,0,1\}$.
		Every Fourier monomial \(e^{\mathrm i m\cdot z}\) occurring in \(F_j\) satisfies \(m_j\ne0\) and \(m_\ell=0\) for all \(\ell>j\). Every Fourier monomial \(e^{\mathrm i n\cdot w}\) occurring in \(G_k\) satisfies \(n_k\ne0\) and \(n_\ell=0\) for all \(\ell<k\).
		\item[(iii)] At the origin,
		\begin{equation}\label{eq:Slarge}
		\sum_{1\le j<k\le N}F_j(0)G_k(0) = (2N-3)/9.
		\end{equation}
	\end{enumerate}
\end{lemma}

The construction below is related to the finite Riesz products used
in \cite[Section~3]{KW} to obtain quantitative versions of Ornstein's
non-inequalities. Here we choose the coefficients so that both
polynomials remain uniformly bounded, while the sum in
\eqref{eq:Slarge} grows linearly with $N$.

\begin{proof}
	We fix the integer $N\geq 2$ and define real numbers
\[
\begin{aligned}
	&s_0:=0,\qquad s_j:=\frac{(-1)^{j+1}}3\quad\text{for all }1\le j\le N,\\
	&r_{N+1}:=0,\qquad r_k:=\frac{(-1)^k}3\quad\text{for all }1\le k\le N.
\end{aligned}
\]
For $z,w\in\T^N$, we set
\begin{align}
	F_j(z)&:=(s_j-s_{j-1})\cos z_j
	\prod_{\ell<j}\frac{1+\cos z_\ell}{2},
	&&\text{for all }1\le j\le N,\label{eq:Fj}\\
	G_k(w)&:=(r_k-r_{k+1})\cos w_k
	\prod_{\ell>k}\frac{1+\cos w_\ell}{2},
	&&\text{for all }1\le k\le N,\label{eq:Gk}\\
	F(z)&:=\sum_{j=1}^N F_j(z),\qquad
	G(w):=\sum_{k=1}^N G_k(w).\label{eq:FG}
\end{align}
By definition, empty products equal one.

	We regard $F$ as a polynomial in the real variables $q_j:=\cos z_j\in[-1,1]$.
	It is affine in each variable separately. Its maximum absolute value on
	the cube is therefore attained at a vertex. If all $q_j=1$, the sum
	telescopes to $s_N$. Otherwise, we denote by $j$ the first index with $q_j=-1$.
	Every subsequent term vanishes, and the value is
	\[
	s_{j-1}-(s_j-s_{j-1})=2s_{j-1}-s_j.
	\]
	For $j=1$ this equals $-1/3$; for $j\ge2$ it equals $1$ or $-1$.
	This proves $\|F\|_\infty\le1$. The identical argument in reverse
	coordinate order proves $\|G\|_\infty\le1$.
	
	The factor $\cos z_j$ in \eqref{eq:Fj} has mean zero, and all remaining
	factors are independent of $z_j$. Thus each $F_j$ has mean zero. Every Fourier monomial $e^{\ii m\cdot z}$ occurring in \(F_j\) has \(m_j\in\{-1,1\}\) and \(m_\ell=0\) for all \(\ell>j\). All remaining coordinate Fourier exponents belong to \(\{-1,0,1\}\). The corresponding assertions for \(G_k\) follow in the same way.
	
	At the origin, $F_j(0)=s_j-s_{j-1}$ and $G_k(0)=r_k-r_{k+1}$. Hence
	\[
	\sum_{j<k}F_j(0)G_k(0)=\sum_{k=2}^N s_{k-1}(r_k-r_{k+1}).
	\]
	Each term with $2\le k<N$ equals $2/9$, and the last term equals $1/9$.
	Thus the sum equals $(2N-3)/9$, proving \eqref{eq:Slarge}.
\end{proof}

%%%%%%%%%%%%%%%%%%%%%%%%%%%

\subsection{An elementary operator-norm lower bound}
For the one-dimensional construction we write $P:=-\ii\,d/dx$; with our
Fourier convention, $(e^{\ii\beta P}\psi)(x)=\psi(x+\beta)$.

For integers $N\ge1$ and $R\ge3$, we define unitaries
\begin{equation}\label{eq:unitaries}
	U_{j,R}:=e^{\ii R^jX},\qquad
	V_{k,R}:=e^{\ii2\pi R^{-k}P},\qquad \text{for all}\ 1\le j,k\le N \,.
\end{equation}
The $U_{j,R}$ commute among themselves, as do the $V_{k,R}$. For $m,n\in\Z^N$, we set $U_R^m:=\prod_jU_{j,R}^{m_j}$ and $V_R^n:=\prod_kV_{k,R}^{n_k}$. For a
trigonometric polynomial $H$ on $\T^{2N}$,
\[
H(z,w)=\sum_{m,n}\widehat H(m,n)e^{\ii(m\cdot z+n\cdot w)},
\]
we define an operator $H_R$ on $L^2(\R)$ by using the following ordering:
\begin{equation}\label{eq:ordered}
	H_R:=\sum_{m,n}\widehat H(m,n)U_R^mV_R^n.
\end{equation}

\begin{lemma}\label{lem:norm}
	For every fixed $N$ and every trigonometric polynomial $H$ on $\T^{2N}$,
	\[
	\liminf_{\substack{R\to\infty\\R\in\mathbb N}}\|H_R\|
	\ge |H(0,0)|.
	\]
\end{lemma}

The strategy in the following proof is sometimes called a F\o lner argument. It uses the familiar fact that translating
a large finite box by one unit changes only a small proportion of its
points; see \cite[Sections~G.3 and~G.5]{BHV}. We place the points of
such a box on the real line and replace them by narrow bumps. Their
locations are chosen so that the multiplication operators also
almost preserve the resulting vector.

\begin{proof}
	We construct unit vectors $\psi_R$ such that
	\begin{equation}\label{eq:approx-eigen}
		\|(U_{j,R}-I)\psi_R\|_2\longrightarrow0,
		\qquad
		\|(V_{k,R}-I)\psi_R\|_2\longrightarrow0
		\quad\text{for all }1\le j,k\le N.
	\end{equation}
	We put $L:=\lfloor\sqrt R\rfloor$ and take the integer $R$ large enough that $L+1<R$.
	We consider the finite set of integers
	\[
	\mathcal D_R:=\left\{\sum_{r=0}^{N-1}d_rR^r:
	d_r\in\{1,\ldots,L\}\right\}
	\]
	with cardinality $L^N$. We fix $\chi\in C_c^\infty((-1,1))$ with
	$\|\chi\|_2=1$, and define
	\begin{equation}\label{eq:test-vector}
		\psi_R(x):=L^{-N/2} R^N \sum_{\ell\in\mathcal D_R}
		\chi\left(R^{2N}\left(x-\frac{2\pi\ell}{R^N}\right)\right).
	\end{equation}
	This is a superposition of $L^N$ bumps with disjoint supports and $\|\psi_R\|_2=1$.
	
	The translation $V_{k,R}$ shifts the centers of the bumps by $-2\pi/R^k$, equivalently
	it replaces $\ell$ by $\ell-R^{N-k}$. The intersection of
	$\mathcal D_R$ with this translate has cardinality
	$(L-1)L^{N-1}$: exactly one base-$R$ digit is shifted by one. Different centers still have disjoint bump supports.
	Consequently
	\begin{equation}\label{eq:translation-error}
		\langle\psi_R,V_{k,R}\psi_R\rangle=1-\frac1L,
		\qquad
		\|(V_{k,R}-I)\psi_R\|_2^2=\frac2L.
	\end{equation}
	
	For a center $x_\ell:=2\pi\ell/R^N$, the number $R^jx_\ell/(2\pi)$
	differs from an integer by
	\[
	\sum_{r=0}^{N-j-1}d_rR^{r-(N-j)}
	\in\left[0,\frac{L}{R-1}\right],
	\]
	with an empty sum if $j=N$. On the corresponding bump, the additional
	phase variation is at most $R^j \times R^{-2N} \le R^{-N}$. Hence
	\begin{equation}\label{eq:modulation-error}
		\|(U_{j,R}-I)\psi_R\|_2
		\le\frac{2\pi L}{R-1}+R^{-N} \longrightarrow0.
	\end{equation}
	Equations \eqref{eq:translation-error} and \eqref{eq:modulation-error}
	prove \eqref{eq:approx-eigen}.
	
	Note that for any unitary $W$, any $m\in\Z$ and any $\psi$ we have
	$\|(W^m-I)\psi\|_2\le |m|\|(W-I)\psi\|_2$. Indeed, this follows from $
	\|(W^{-1}-I)\psi\|=\|(W-I)\psi\|$ by telescoping.
	
	Telescoping the finite product and using \eqref{eq:approx-eigen}, we obtain
	\[
	\|(U_R^mV_R^n-I)\psi_R\|_2\longrightarrow0
	\qquad\text{for every fixed }m,n.
	\]
	Since the sum in \eqref{eq:ordered} is finite, it follows that
	\[
	\|H_R\psi_R-H(0,0)\psi_R\|_2\longrightarrow0.
	\]
	Taking norms proves the lemma.
\end{proof}

%%%%%%%%%%%%%%%%%%%%%%%%%

\subsection{Proof of Theorem~\ref{thm:counterexample}}

Finally, we can prove the main result of this section.

\begin{proof}
	\textit{Step 1.} We begin with the case of dimension $d=1$ and construct, for any $M>0$, real smooth bounded periodic functions $a,b$ on $\R$ with
	\[
	\|a'\|_\infty\leq1,\qquad\|b'\|_\infty\leq1,
	\qquad\|[a(X),b(P)]\|>M \,.
	\]
		
	We choose the polynomials $F,G$ from Lemma~\ref{lem:polynomials} for a given integer $N\geq 2$. By property (ii), $F,G$ have mean zero. We write
	\[
	F(z)=\sum_m\widehat F(m)e^{\ii m\cdot z},\qquad
	G(w)=\sum_n\widehat G(n)e^{\ii n\cdot w},
	\]
	and, for all \(m,n\) with \(\widehat F(m)\ne0\) and \(\widehat G(n)\ne0\), we put
	\[
	\alpha_m(R):=\sum_{j=1}^N m_jR^j,
	\qquad
	\beta_n(R):=2\pi\sum_{k=1}^N n_kR^{-k}.
	\]
	These numbers are nonzero for $R\ge3$: the highest nonzero term in
	$\sum_jm_jR^j$ dominates the sum of all lower terms, and reversing the
	indices gives the analogous assertion for $\sum_kn_kR^{-k}$.
	We define
	\begin{equation}\label{eq:symbols}
		a_R(x):=\sum_{\widehat F(m)\ne0}\frac{\widehat F(m)}{\ii\alpha_m(R)}
		e^{\ii\alpha_m(R)x},
		\qquad
		b_R(\xi):=\sum_{\widehat G(n)\ne0}\frac{\widehat G(n)}{\ii\beta_n(R)}
		e^{\ii\beta_n(R)\xi}.
	\end{equation}
	Reality of $F,G$ implies reality of $a_R,b_R$, by pairing the terms corresponding to \(m\) and \(-m\), and to \(n\) and \(-n\). The functions $a_R$ and $b_R$ are finite trigonometric polynomials, hence smooth and bounded; $2\pi/R$ and $R^N$ are periods, respectively. Furthermore,
	\begin{align*}
		a_R'(x)&=F(Rx,R^2x,\ldots,R^Nx),\\
		b_R'(\xi)&=G(2\pi R^{-1}\xi,2\pi R^{-2}\xi,\ldots,2\pi R^{-N}\xi),
	\end{align*}
	so, by (i) in Lemma \ref{lem:polynomials},
	\begin{equation}\label{eq:gradient}
		\|a_R'\|_\infty\le1,\qquad\|b_R'\|_\infty\le1.
	\end{equation}
	
	We set $h(t):=(e^{\ii t}-1)/(\ii t)$ for $t\ne0$, and $h(0):=1$.
	Since
	\[
	[e^{\ii\alpha X},e^{\ii\beta P}]
	=(1-e^{\ii\alpha\beta})e^{\ii\alpha X}e^{\ii\beta P},
	\]
	the division by $\ii\alpha$ and $\ii\beta$ in \eqref{eq:symbols} gives
	the identity
	\begin{equation}\label{eq:commutator}
		-\ii[a_R(X),b_R(P)]
		=\sum_{m,n}\widehat F(m)\widehat G(n)
		h\bigl(\alpha_m(R)\beta_n(R)\bigr)U_R^mV_R^n.
	\end{equation}
	Here, strictly speaking, the sum ranges over \(m\) and \(n\) with $\widehat F(m)\neq 0$ and $\widehat G(n)\neq 0$ (because $\alpha_m(R)$ and $\beta_n(R)$ are only defined for such $m$ and $n$).
	
	Let us consider $j$, $k$, $m$ and $n$ such that \(\widehat F_j(m)\ne0\) and \(\widehat G_k(n)\ne0\). Property (ii) in Lemma \ref{lem:polynomials} gives
	\[
	\alpha_m(R)=m_jR^j(1+O(R^{-1})),\qquad
	\beta_n(R)=2\pi n_kR^{-k}(1+O(R^{-1})).
	\]
	Here $m_j,n_k\in\{-1,1\}$. If $j<k$, the product $\alpha_m(R)\beta_n(R)$ tends to zero as $R\to\infty$ and $h(\alpha_m(R)\beta_n(R))\to1$. If $j>k$, its absolute value tends to
	infinity and $h(\alpha_m(R)\beta_n(R))\to0$, since $|h(t)|\le2/|t|$.
	If $j=k$, the product tends to $2\pi m_jn_j=\pm2\pi$, at which $h$
	vanishes. 
	
	Consequently, if we decompose $F$ and $G$ as in Lemma \ref{lem:polynomials}, then only the terms with $j<k$ survive in the limit $R\to\infty$. Using
	\[
	\sum_{m,n} \sum_{j<k} \widehat F_j(m)\widehat G_k(n)
	U_R^mV_R^n = S_R
	\]
	with the  trigonometric polynomial $S$ on $\T^{2N}$, defined by
	\begin{equation}\label{eq:S}
		S(z,w):=\sum_{1\le j<k\le N}F_j(z)G_k(w),
	\end{equation}
	we arrive at	
	\begin{equation}\label{eq:triangular}
		\left\|-\ii[a_R(X),b_R(P)]-
		S_R \right\|\longrightarrow0.
	\end{equation}
	This is operator-norm convergence: for fixed \(N\), there are only finitely many Fourier terms, their scalar coefficients converge, and \(\|U_R^mV_R^n\|=1\) because all factors are unitary.
	
	Lemma~\ref{lem:norm} and \eqref{eq:Slarge} now prove
	\begin{equation}\label{eq:divergence}
		\liminf_{\substack{R\to\infty\\R\in\mathbb N}}
		\|[a_R(X),b_R(P)]\|
		\ge S(0,0)= (2N-3)/9.
	\end{equation}
	We now choose $N$ so large that $(2N-3)/9>M$. For all sufficiently
	large integers $R$, \eqref{eq:gradient} and \eqref{eq:divergence}
	then give the required functions $a:=a_R$ and $b:=b_R$.

\medskip
	
	\textit{Step 2.} 
	To obtain the conclusion in general dimension $d\geq 1$, we put $A(x):=a(x_1)$ and
	$B(\xi):=b(\xi_1)$ on $\R^d$. Their commutator is the one-dimensional
	commutator tensored with the identity, so its norm is unchanged, and
	$\|\nabla A\|_\infty,\|\nabla B\|_\infty\leq1$.
	We choose a real $\chi\in C_c^\infty(\R^d)$ with $0\leq\chi\leq1$ and
	$\chi=1$ near zero, and set
	\[
	A_L(x):=\chi(x/L)A(x),\qquad B_L(\xi):=\chi(\xi/L)B(\xi).
	\]
	These are real smooth compactly supported functions, with
	\[
	\begin{split}
		\|\nabla A_L\|_\infty&\leq1+L^{-1}\|\nabla\chi\|_\infty\|A\|_\infty,\\
		\|\nabla B_L\|_\infty&\leq1+L^{-1}\|\nabla\chi\|_\infty\|B\|_\infty.
	\end{split}
	\]
	The uniformly bounded operators $A_L(X)$ and $B_L(-\ii\nabla)$ converge
	strongly to $A(X)$ and $B(-\ii\nabla)$, respectively. Their commutators
	therefore converge strongly, and
	\[
	\liminf_{L\to\infty}\|[A_L(X),B_L(-\ii\nabla)]\|
	\geq\|[A(X),B(-\ii\nabla)]\|>M.
	\]
	We divide $A_L,B_L$ by the respective upper bounds for their gradient norms.
	The divisors tend to one, so for sufficiently large $L$ the normalized
	functions have gradient norms at most one and commutator norm greater than
	$M$. Since $M$ is arbitrary, this disproves the uniform estimate in
	Theorem~\ref{thm:counterexample}.
\end{proof}

The use of widely separated frequencies in Step 1 of the preceding proof is related to Bourgain's
argument in \cite[proof of Lemma~1]{Bourgain}, where the largest
frequency determines the action of the Hilbert transform. Here we
arrange the position and momentum frequencies in opposite orders.
This makes the commutator retain precisely the terms with $j<k$
in the limit.

\begin{remark}
	For the periodic functions $a,b$ constructed above, we put
	$A(x):=(x+a(x))/2$ and $B(\xi):=(\xi+b(\xi))/2$. Then $0\leq A',B'\leq1$
	and, as a form on pairs of Schwartz functions,
	\[
	-\ii[A(X),B(P)]
	=\frac14\bigl(-\ii[a(X),b(P)]+a'(X)+b'(P)+I\bigr).
	\]
	The right side is bounded and has norm at least
	$(\|[a(X),b(P)]\|-3)/4$. Thus the commutator norm $\|[A(X),B(P)]\|$ is not uniformly bounded even among smooth functions $A,B$ with $0\leq A',B'\leq1$.
\end{remark}

\begin{remark}\label{rem:unbounded}
	For every $d\geq1$, Theorem~\ref{thm:counterexample} also implies that
	there are real globally Lip\-schitz functions $a,b$ on $\R^d$ for which
	the commutator form has no bounded extension to $L^2(\R^d)$.
	Indeed, since adding constants does not affect the commutator, we may work in the real Banach space of Lipschitz functions $f$ on \(\mathbb R^d\) vanishing at the origin, equipped with the norm \(\|\nabla f\|_\infty\). The commutator form is well defined on pairs of Schwartz functions and depends continuously on each function in this norm, since
	\[
	|a(x)|\leq\|\nabla a\|_\infty|x|,
	\qquad |b(\xi)|\leq\|\nabla b\|_\infty|\xi|.
	\]
	If every pair had a bounded commutator extension, the closed graph
	theorem would make the map from the two functions to that extension
	separately bounded: convergence of one function in the Lipschitz norm
	implies convergence of its commutator forms on every pair of Schwartz
	functions, and therefore identifies any limit in operator norm.
	The uniform boundedness principle would then give a uniform estimate
	by a constant times $\|\nabla a\|_\infty\|\nabla b\|_\infty$,
	contrary to Theorem~\ref{thm:counterexample}.
\end{remark}

%%%%%%%%%%%%%%%%%%%%%%%%%%%%%%%%
%%%%%%%%%%%%%%%%%%%%%%%%%%%%%%%%

\section{Failure of the estimates at and below $p=1$}\label{sec:counter}

We prove Theorem~\ref{thm:counter}. One-dimensional constructions give
large trace norms and settle the weak endpoint $p=1$. Applying orthogonal
projections in the remaining variables yields the counterexamples in
every dimension. We also give individual nonmembership examples and
explain the distinction between dimension one and higher dimensions.

\subsection{The one-dimensional construction}\label{sec:counter-longitudinal}

We first construct one-dimensional operators whose trace norms grow
with the number of frequencies. All constants below are independent of
the truncation parameters.

We choose a nonzero nonnegative real $\rho\in C_c^\infty(\R)$, set $\widetilde\rho(t):=\rho(-t)$ and put $\chi:=\rho*\widetilde\rho$.
Then $\widehat\chi\geq0$ and $\widehat\chi$ is strictly positive near zero.
We set $N_j:=8^j$ and
\begin{equation}\label{eq:countera}
 A_L(t):=\chi(t)\sum_{j=1}^L\frac{\cos(N_jt)}{jN_j},
 \qquad
 A(t):=\chi(t)\sum_{j=1}^\infty\frac{\cos(N_jt)}{jN_j}.
\end{equation}
The series converges uniformly. Moreover,
\begin{equation}\label{eq:counter-w12}
 \sup_L\bigl(\|A_L\|_\infty+\|A_L'\|_2\bigr)<\infty,
 \qquad A\in C_c(\R)\cap W^{1,2}(\R).
\end{equation}
Indeed, the differentiated series with $\chi'$ converges absolutely and
uniformly. The other series is
$-\chi(t)\sum_j\sin(N_jt)/j$; it converges in $L^2$ because the sines
are orthogonal on every interval of length $2\pi$, and the fixed support
of $\chi$ is covered by finitely many such intervals.
Consequently, for every $p>0$,
\begin{equation}\label{eq:countergrad}
 \sup_L\bigl(\|A_L\|_{\ell^p(L^2)}+
                  \|A_L'\|_{\ell^p(L^2)}\bigr)<\infty.
\end{equation}

We denote by $P_I$ the orthogonal projection onto Fourier frequencies in
an interval $I$. For $L\geq J$, including $L=\infty$ with $A_\infty:=A$,
we put
\[
 H_{J,L}:=P_{(1,N_J)}A_L(X)P_{(-N_J,0)}.
\]
This operator is trace class: we factor the multiplication by $A_L$ into
two factors $|A_L|^{1/2}$ and a bounded sign factor. Each outer factor,
paired with its frequency projection, is Hilbert--Schmidt because
$A_L\in L^1$ and both frequency intervals have finite measure.

\begin{lemma}\label{lem:counter-longitudinal}
For $L\geq J$ one has
\begin{equation}\label{eq:counter-trace}
 \|H_{J,L}\|_{\Sp^1}\gtrsim\sum_{j=1}^J\frac1j.
\end{equation}
Consequently, for $0<p\leq1$,
\[
 \|H_{J,L}\|_{\Sp^p}\gtrsim\sum_{j=1}^J\frac1j,
\]
and, for $0<p<1$,
\begin{equation}\label{eq:counter-weak}
 \|H_{J,L}\|_{\Sp^p_{\rm weak}}
       \gtrsim_p\sum_{j=1}^J\frac1j.
\end{equation}
\end{lemma}

The operators $H_{J,L}$ are obtained from integral Hankel operators
by restricting the input and output frequencies. The distinction between square summability of lacunary
derivative coefficients and the summability required for trace class
is classical in this setting; see \cite{Peller80} and
\cite[Sections~4.1 and~5.3]{Peller24}. We give a direct proof
that retains the frequency restrictions needed for the
higher-dimensional construction.

\begin{proof}
We put
\[
 I_j:=(-N_j/2,-N_j/4),\qquad
 V_J:=\sum_{j=1}^J e^{\ii N_jX}P_{I_j}.
\]
Each summand shifts Fourier support from $I_j$ to $I_j+N_j$.
Both families of intervals are pairwise disjoint and lie in the
respective input and output bands of $H_{J,L}$, so $\|V_J\|=1$.
For a bounded interval $I$ and a function $F\in L^1(\R)\cap L^\infty(\R)$,
the Hilbert--Schmidt factorization used above gives
\[
 \Tr\bigl(P_I F(X)P_I\bigr)
 =(2\pi)^{-1}|I|\int_{\R}F(x)\,dx
 =(2\pi)^{-1/2}|I|\widehat F(0).
\]
Moreover, $\widehat\chi\geq0$ and the $j$-th summand in
\eqref{eq:countera} imply
\[
 \widehat A_L(N_j)\geq\frac{\widehat\chi(0)}{2jN_j}
 \qquad\text{for all }j\leq J.
\]
This also holds for $L=\infty$ by uniform convergence and fixed
compact support. Each summand of $V_J^*H_{J,L}$ has range in
$P_{I_j}L^2$, so inserting $P_{I_j}$ on its right does not change its
trace. Using $|I_j|=N_j/4$, we obtain
\begin{align*}
 \|H_{J,L}\|_{\Sp^1}
 &\geq |\Tr(V_J^*H_{J,L})|\\
 &=\left|\sum_{j=1}^J
 \Tr\bigl(P_{I_j}(e^{-\ii N_j\cdot}A_L)(X)P_{I_j}\bigr)\right|\\
 &=(2\pi)^{-1/2}\sum_{j=1}^J|I_j|\widehat A_L(N_j)
 \gtrsim\sum_{j=1}^J\frac1j,
\end{align*}
which proves \eqref{eq:counter-trace}.
For $p\leq1$, $\|T\|_{\Sp^1}\leq\|T\|_{\Sp^p}$, and for $p<1$,
\[
 \|T\|_{\Sp^1}
 \leq\left(\sum_{\nu\geq1}\nu^{-1/p}\right)
       \|T\|_{\Sp^p_{\rm weak}},
\]
which proves the remaining assertions.
\end{proof}

%%%%%%%%%%%%%%%%%%%%%%%%%

\subsection{Projections in the remaining variables}\label{sec:counter-transverse}

For the higher-dimensional construction, we fix $d\geq2$ and write
$x=(x_1,x')$ and $\xi=(\xi_1,\xi')$. We choose
$f\in C_c^\infty(\R^{d-1})$ with $f\geq0$ and $f\geq1$ on
$[-\pi,\pi]^{d-1}$.

\begin{lemma}\label{lem:counter-transverse}
For every $R\geq8$ there is a finite-dimensional subspace
$E_R\subset L^2(\R^{d-1})$, all of whose functions have Fourier support
in $[-R,R]^{d-1}$, such that
\[
 \dim E_R\asymp_d R^{d-1},
 \qquad P_{E_R}f(X')P_{E_R}\geq c I_{E_R},
\]
where $c>0$ is independent of $R$.
\end{lemma}

The following construction is an elementary instance of the
concentration phenomenon studied by Landau and Pollak
\cite{LandauPollak}. We construct a space of dimension comparable
to $R^{d-1}$ whose functions have Fourier support in a prescribed
cube and retain a fixed proportion of their $L^2$ mass in a
fixed spatial cube.

\begin{proof}
We fix $0<\delta<1/4$, and define
\[
 \widehat\theta_0:=(2\delta)^{-(d-1)/2}\1_{[-\delta,\delta]^{d-1}},
 \qquad \theta_k(x'):=e^{\ii k\cdot x'}\theta_0(x'),
\]
for all $k\in\Z^{d-1}\cap[-R/2,R/2]^{d-1}$.
These functions are orthonormal, since their Fourier supports are
disjoint. They span a space of dimension comparable to $R^{d-1}$, with the
required Fourier support. The inverse Fourier transform of the small
cube indicator satisfies $|\theta_0(x')|\geq c_0>0$ on
$[-\pi,\pi]^{d-1}$: each factor is a constant times
$\sin(\delta x_j)/x_j$, with its continuous value at zero.
Thus, for any coefficients $z_k$,
\begin{align*}
 \int_{\R^{d-1}} f(x')\left|\sum_kz_k\theta_k(x')\right|^2\,dx'
 &\geq c_0^2\int_{[-\pi,\pi]^{d-1}}
                  \left|\sum_kz_ke^{\ii k\cdot x'}\right|^2\,dx'\\
 &=c_0^2(2\pi)^{d-1}\sum_k|z_k|^2.
\end{align*}
This proves the operator inequality.
\end{proof}

We write $D_R:=P_{E_R}f(X')P_{E_R}$, regarded as an operator on $E_R$,
and $M_R:=\dim E_R$. For any compact $H$ and $p>0$, the tensor-product
singular values give
\begin{equation}\label{eq:counter-tensor-strong}
 \|H\otimes D_R\|_{\Sp^p}
 \geq c M_R^{1/p}\|H\|_{\Sp^p}.
\end{equation}
They also give
\begin{equation}\label{eq:counter-tensor-weak}
 \|H\otimes D_R\|_{\Sp^p_{\rm weak}}
 \geq c M_R^{1/p}\|H\|_{\Sp^p_{\rm weak}}.
\end{equation}
Indeed, at least $jM_R$ singular values of $H\otimes D_R$ are at least
$c s_j(H)$.

We choose a fixed smooth nondecreasing function $h$ with $h=0$ on
$(-\infty,0]$ and $h=1$ on $[1,\infty)$. We define
\begin{equation}\label{eq:counter-box}
 g_R(\xi'):=\prod_{\nu=2}^d
       h(\xi_\nu+R+1)h(R+1-\xi_\nu),
 \qquad
 b_R(\xi):=h(\xi_1)h(2R-\xi_1)g_R(\xi').
\end{equation}
Thus $g_R=1$ on $[-R,R]^{d-1}$ and vanishes outside
$[-R-1,R+1]^{d-1}$. All transition layers have width one.

Every component of $\nabla b_R$ is uniformly bounded and supported in
boundary layers meeting at most $O_d(R^{d-1})$ unit cubes. The two
longitudinal end layers meet $O(R^{d-1})$ cubes. Each transverse layer
meets $O(R\,R^{d-2})=O(R^{d-1})$ cubes. Consequently,
\begin{equation}\label{eq:counter-box-grad}
 \|\nabla b_R\|_{\ell^p(L^2)}^p\lesssim_{d,p}R^{d-1}
 \qquad\text{for all }0<p<\infty.
\end{equation}
The estimate holds uniformly after translation as well, since a unit
cube in one translated grid meets a bounded number of cubes in the
original grid.

We will use the following identity involving these projections. For bounded measurable functions
$\alpha,\beta$ on $\R$, we set
$a(x):=\alpha(x_1)f(x')$ and
$b(\xi):=\beta(\xi_1)g_R(\xi')$. Then, on
$L^2(\R)\otimes E_R$,
\begin{equation}\label{eq:counter-product-compression}
 \begin{split}
 &(I\otimes P_{E_R})[a(X),b(-\ii\nabla)](I\otimes P_{E_R}) =[\alpha(X_1),\beta(-\ii\partial_1)]\otimes D_R.
 \end{split}
\end{equation}
Indeed, the Fourier-support property of $E_R$ gives
$g_R(-\ii\nabla')P_{E_R}=P_{E_R}=P_{E_R}g_R(-\ii\nabla')$,
and the identity follows by expanding the commutator.

In particular, for $R=N_J$, $L\geq J$, and
$a_L(x):=A_L(x_1)f(x')$, the function $h(t)h(2R-t)$ is zero on
$(-R,0)$ and one on $(1,R)$. Thus
\eqref{eq:counter-product-compression} gives
\begin{equation}\label{eq:counter-compression}
 (P_{(1,R)}\otimes P_{E_R})[a_L(X),b_R(-\ii\nabla)]
 (P_{(-R,0)}\otimes P_{E_R})=-H_{J,L}\otimes D_R.
\end{equation}

%%%%%%%%%%%%%%%%%%%%%%%%%%%%

\subsection{Proof of Theorem~\ref{thm:counter}}
Finally, we can prove the main result of this section.

\begin{proof}[Proof of Theorem~\ref{thm:counter}]
\emph{Failure of the strong estimate for $0<p\leq1$ and the weak estimate for $0<p<1$.}

We first suppose that $d\geq2$. We take $R=N_J$, $a_J=A_J\otimes f$, and $b_R$
as above. These are real smooth compactly supported functions. The norm of \(\nabla a_J\) stays bounded, by \eqref{eq:counter-w12} and fixed compact
support, while \eqref{eq:counter-box-grad} gives
\[
 \|\nabla a_J\|_{\ell^p(L^2)}\lesssim_{d,p}1,
 \qquad \|\nabla b_R\|_{\ell^p(L^2)}
           \lesssim_{d,p}R^{(d-1)/p}.
\]
Schatten norms and weak Schatten quasinorms cannot increase under
multiplication on either side by an orthogonal projection. Applying
\eqref{eq:counter-compression}, Lemma~\ref{lem:counter-longitudinal}, and
\eqref{eq:counter-tensor-strong} yields
\[
 \|[a_J(X),b_R(-\ii\nabla)]\|_{\Sp^p}
 \gtrsim_{d,p}R^{(d-1)/p}\sum_{j=1}^J\frac1j,
 \qquad\text{for all }0<p\leq1.
\]
The ratio to the product of the gradient norms therefore tends to
infinity, disproving the strong estimate. In dimension one the same proof omits $f,g_R,E_R,D_R$;
the derivative of $h(\xi)h(2N_J-\xi)$ consists of two fixed-width
transition layers and has uniformly bounded $\ell^p(L^2)$ norm.

For $0<p<1$, \eqref{eq:counter-weak} and
\eqref{eq:counter-tensor-weak} give the same lower bound for the
weak Schatten quasinorm. Since
$\|\nabla b_R\|_{\ell^p_{\rm weak}(L^2)}
\leq\|\nabla b_R\|_{\ell^p(L^2)}$, the weak estimate in
Theorem~\ref{thm:counter} fails as well.

\medskip

\emph{Failure of the weak estimate at $p=1$.}
We use the following elementary inequality for compact operators:
\begin{equation}\label{eq:weaktrace-hs}
 \|T\|_{\Sp^2}^2\leq2\|T\|\,\|T\|_{\Sp^1_{\rm weak}}.
\end{equation}
Indeed, with $n_T(t):=\#\{j:s_j(T)>t\}$, we have
$n_T(t)\leq t^{-1}\|T\|_{\Sp^1_{\rm weak}}$ and
\[
 \|T\|_{\Sp^2}^2
 =\int_0^{\|T\|}2t\,n_T(t)\,dt
 \leq2\|T\|\,\|T\|_{\Sp^1_{\rm weak}}.
\]

We first give the construction in dimension one. We choose an even
$\kappa\in C_c^\infty((-2,2))$ with $0\leq\kappa\leq1$ and
$\kappa=1$ on $[-1,1]$. For integers $N\geq8$, we set
\begin{equation}\label{eq:weaktrace-symbol}
 L_N:=\operatorname{arsinh}N,\qquad
 \beta_N(t):=\operatorname{arsinh}t\,
       \kappa\left(\frac{\operatorname{arsinh}t}{L_N}\right).
\end{equation}
Then $\beta_N\in C_c^\infty(\R)$ is real and
\begin{equation}\label{eq:weaktrace-symbol-bounds}
 \begin{split}
 \beta_N(t)&=\operatorname{arsinh}t\quad\text{for all }|t|\leq N,\\
 \|\beta_N\|_\infty&\leq2L_N,\qquad
 \supp\beta_N\subset[-3N^2,3N^2],\qquad
 |\beta_N'(t)|\lesssim(1+|t|)^{-1}.
 \end{split}
\end{equation}
For the support assertion, we use
$\sinh(2L_N)=2N\sqrt{1+N^2}\leq3N^2$.
The derivative estimate follows by the chain rule; its constant depends
only on $\kappa$. In particular,
\begin{equation}\label{eq:weaktrace-beta-gradient}
 \|\beta_N'\|_{\ell^1_{\rm weak}(L^2)}\lesssim1,
\end{equation}
since its local $L^2$ norm on $[k,k+1)$ is bounded by a constant
times $(1+|k|)^{-1}$.

With $\chi$ as in Subsection~\ref{sec:counter-longitudinal}, we put
\[
 a_N(x):=N^{-1}\chi(x)\cos(Nx),\qquad b_N:=\beta_N,
 \qquad T_N:=[a_N(X),b_N(-\ii\nabla)].
\]
The derivative of $a_N$ has fixed compact support and uniformly bounded
$L^2$ norm, so $\|a_N'\|_{\ell^1(L^2)}\lesssim1$. Also,
\begin{equation}\label{eq:weaktrace-op}
 \|T_N\|\leq2\|a_N\|_\infty\|b_N\|_\infty
             \lesssim L_N/N.
\end{equation}
The Fourier kernel of the commutator gives
\begin{equation}\label{eq:weaktrace-hs-kernel}
 \|T_N\|_{\Sp^2}^2
 =(2\pi)^{-1}\iint_{\R^2}|\widehat a_N(z)|^2
       |b_N(\eta+z)-b_N(\eta)|^2\,d\eta\,dz.
\end{equation}
We choose a fixed sufficiently small $\delta>0$ and restrict the integral
to $|z-N|<\delta$ and $-N/2<\eta<-N/3$. Fourier positivity gives
\[
 \widehat a_N(z)
 =\frac{\widehat\chi(z-N)+\widehat\chi(z+N)}{2N}
 \gtrsim N^{-1}.
\]
Both $|\eta|$ and $\eta+z$ are comparable to $N$ and less than $N$.
Thus \eqref{eq:weaktrace-symbol-bounds} implies
\[
 |b_N(\eta+z)-b_N(\eta)|
 =\operatorname{arsinh}(\eta+z)+\operatorname{arsinh}(-\eta)
 \gtrsim L_N.
\]
The region of integration has measure comparable to $N$, whence
\[
 \|T_N\|_{\Sp^2}^2\gtrsim L_N^2/N.
\]
Combining this with \eqref{eq:weaktrace-hs} and
\eqref{eq:weaktrace-op} yields
\begin{equation}\label{eq:weaktrace-divergence}
 \|T_N\|_{\Sp^1_{\rm weak}}\gtrsim L_N\asymp\log N.
\end{equation}
The two gradient norms on the right side of the proposed estimate
stay bounded, by \eqref{eq:weaktrace-beta-gradient} and the bound on
$a_N'$. This proves its failure in dimension one.

For $d\geq2$, we use the functions $f,g_R$ from Subsection~\ref{sec:counter-transverse} and set
\[
 R:=N^3,\qquad
 a_N(x):=N^{-1}\chi(x_1)\cos(Nx_1)f(x'),\qquad
 b_N(\xi):=\beta_N(\xi_1)g_R(\xi').
\]
The gradient of $a_N$ has fixed compact support and uniformly bounded
$L^2$ norm, so $\|\nabla a_N\|_{\ell^1(L^2)}\lesssim_d1$.
The local $L^2$ norms of $\partial_1b_N$ are bounded by those of
$\beta_N'$, repeated on $O_d(R^{d-1})$ transverse cubes. Thus
\eqref{eq:weaktrace-beta-gradient} gives
\[
 \|\partial_1b_N\|_{\ell^1_{\rm weak}(L^2)}\lesssim_d R^{d-1}.
\]
Each transverse derivative has amplitude $O(L_N)$ and is supported in
boundary layers meeting $O_d(N^2R^{d-2})$ unit cubes. The quasitriangle
inequality for the weak sequence norm therefore gives
\begin{equation}\label{eq:weaktrace-gradient-d}
 \|\nabla b_N\|_{\ell^1_{\rm weak}(L^2)}
 \lesssim_d R^{d-1}+N^2L_NR^{d-2}\lesssim_d R^{d-1}.
\end{equation}
Here we used $R=N^3$ and $L_N\lesssim N$.
Identity \eqref{eq:counter-product-compression},
the tensor-product bound \eqref{eq:counter-tensor-weak}, and
\eqref{eq:weaktrace-divergence} imply
\[
 \|[a_N(X),b_N(-\ii\nabla)]\|_{\Sp^1_{\rm weak}}
 \gtrsim_d R^{d-1}L_N.
\]
The ratio to the product of the gradient norms tends to infinity.
This proves the remaining assertion of Theorem~\ref{thm:counter}.
\end{proof}

\begin{remark}\label{rem:counter-one-dimensional}
	For every $d\geq1$ and $0<p<1$, our proof shows that even the weaker estimate
	\[
	\|[a(X),b(-\ii\nabla)]\|_{\Sp^p_{\rm weak}}
	\lesssim_{d,p}\|\nabla a\|_{\ell^p(L^2)}
	\|\nabla b\|_{\ell^p(L^2)}
	\]
	fails for real $a,b\in C_c^\infty(\R^d)$.
	Indeed, for $d\geq2$, the functions in the first part of the proof satisfy,
	with $R=N_J$,
	\[
	\|[a_J(X),b_R(-\ii\nabla)]\|_{\Sp^p_{\rm weak}}
	\gtrsim_{d,p}R^{(d-1)/p}\sum_{j=1}^J\frac1j,
	\qquad
	\|\nabla a_J\|_{\ell^p(L^2)}
	\|\nabla b_R\|_{\ell^p(L^2)}
	\lesssim_{d,p}R^{(d-1)/p}.
	\]
	Their ratio tends to infinity. In dimension one, the same argument
	applies with the factor $R^{(d-1)/p}$ replaced by one.
\end{remark}

%%%%%%%%%%%%%%%%%%%%%

\subsection{Individual nonmembership examples}\label{sec:counter-individual}

The constructions in this section also give the following stronger conclusions about
membership for fixed functions.

\begin{proposition}\label{prop:counter-individual}
	For every $d\geq1$ and $0<p\leq1$, there are bounded real functions $a,b$ on
	$\R^d$ whose distributional gradients satisfy
	\[
	\nabla a,\nabla b\in\ell^p(L^2),
	\qquad [a(X),b(-\ii\nabla)]\in\Sp^2\setminus\Sp^p.
	\]
	If $0<p<1$, the functions can be chosen so that the commutator does
	not belong to $\Sp^p_{\rm weak}$.
	In dimension one, the same pair $a,b$ can be used for every
	$0<p\leq1$, with both derivatives compactly supported and in $L^2$.
\end{proposition}

\begin{proof}
Since $\ell^p(L^2)\subset L^2$ for $p\leq1$, the gradient assumptions
already imply that the commutator belongs to $\Sp^2$, by part~(b) of
Theorem~\ref{thm:main} with exponent $2$.

\emph{Individual counterexamples in dimensions $d\geq2$.}

We fix $d\geq2$ and $0<p\leq1$. We use the fixed function
\[
 a(x):=A(x_1)f(x')
\]
with the infinite series in \eqref{eq:countera}.
Its gradient has compact support and belongs to $L^2$, hence to every
$\ell^q(L^2)$, $q>0$.

We choose strictly increasing integers $J_k$ so large that
\begin{equation}\label{eq:counter-choice}
 \sum_{j=1}^{J_k}\frac1j\geq k2^k,
 \qquad R_k:=N_{J_k},\qquad
 \alpha_k:=2^{-k}R_k^{-(d-1)/p}.
\end{equation}
We choose positive integers $q_k\to\infty$ so that the intervals
\[
 [q_k-R_k-2,q_k+2R_k+2]
\]
are pairwise disjoint and separated by positive gaps. This can be done
inductively, always putting the next interval to the right of the
preceding ones. We set
\begin{equation}\label{eq:counter-glued-b}
 b(\xi):=\sum_{k\geq1}\alpha_k b_{R_k}(\xi_1-q_k,\xi').
\end{equation}
The supports in the first coordinate are pairwise disjoint and form a locally finite family. Thus $b$ is real and smooth; since $\alpha_k\to0$, it is bounded and
vanishes at infinity. The boundary-layer count gives
\begin{equation}\label{eq:counter-glued-grad}
 \|\nabla b\|_{\ell^p(L^2)}^p
 \lesssim_{d,p}\sum_k\alpha_k^pR_k^{d-1}
 =\sum_k2^{-kp}<\infty.
\end{equation}
This also follows from the triangle inequality on each cube and
$(\sum t_k)^p\leq\sum t_k^p$ for nonnegative $t_k$.

On the input frequencies
$(q_k-R_k,q_k)\times[-R_k,R_k]^{d-1}$ the function $b$ is zero.
On the output frequencies
$(q_k+1,q_k+R_k)\times[-R_k,R_k]^{d-1}$ it equals $\alpha_k$.
All other summands vanish on both rectangles by the separation of
the enlarged intervals. Modulation by $e^{\ii q_kx_1}$ commutes with
$A(X_1)$, so multiplying the commutator on the left and right by the
corresponding Fourier projections and by $I\otimes P_{E_{R_k}}$ gives an
operator unitarily equivalent to
\[
 -\alpha_k H_{J_k,\infty}\otimes D_{R_k}.
\]
Its $\Sp^p$ norm is at least
\[
 c_{d,p}\alpha_k R_k^{(d-1)/p}
                    \sum_{j=1}^{J_k}\frac1j
 \geq c_{d,p}k.
\]
The full commutator therefore cannot belong to $\Sp^p$. For $p<1$,
the same argument with the weak norm proves nonmembership in
$\Sp^p_{\rm weak}$. We have used only the fact that multiplication on either
side by an orthogonal projection cannot increase either quasinorm.

\medskip

\emph{The one-dimensional example.}
We set $a:=A$ and $b:=h$, with the functions defined above. Then $a$ is
bounded and compactly supported, $a'\in L^2$ has compact support, and
$b'\in C_c^\infty(\R)$. For every $J$,
\[
 P_{(1,N_J)}[A(X),h(-\ii\nabla)]P_{(-N_J,0)}=-H_{J,\infty}.
\]
By Lemma~\ref{lem:counter-longitudinal}, the trace norms of these operators tend
to infinity. Thus the commutator is not trace class.
Both derivatives belong to $\ell^p(L^2)$ for every $p>0$, whereas the
commutator belongs to no $\Sp^p$ with $p\leq1$, and to no
$\Sp^p_{\rm weak}$ with $p<1$, by the inclusions used in
Lemma~\ref{lem:counter-longitudinal}. This proves the remaining case of the proposition.
\end{proof}

\begin{remark}\label{rem:counter-compact-gradients}
The assertion about compactly supported derivatives in
Proposition~\ref{prop:counter-individual} is specific to dimension one.
Indeed, if $d\geq2$ and $a,b\in W^{1,2}_{\loc}(\R^d)$ have compactly
supported gradients, then each function is constant on the connected
exterior of a sufficiently large ball. Thus there are constants $c_a,c_b$
such that $a-c_a$ and $b-c_b$ are compactly supported and belong to $L^2$,
and hence to $\ell^q(L^2)$ for every $q>0$. The Birman--Solomyak estimate
\eqref{eq:bs}, applied to these functions, and the inclusions between
Schatten classes therefore give
\[
 [a(X),b(-\ii\nabla)]
 =[(a-c_a)(X),(b-c_b)(-\ii\nabla)]\in\Sp^q
 \qquad\text{for every }q>0.
\]
\end{remark}

%%%%%%%%%%%%%%%%%%%%%%%%%%%%%%%%%%%%%%
%%%%%%%%%%%%%%%%%%%%%%%%%%%%%%%%%%%%%%

\appendix
\setcounter{theorem}{0}
\renewcommand{\thetheorem}{A.\arabic{theorem}}
\section{Poincar\'e inequalities and approximation}
\label{app:domains}

All cubes in this appendix have sides parallel to the coordinate axes.
For a cube $Q$ we write $\ell(Q)$ for its side length and
$f_Q:=|Q|^{-1}\int_Q f\,dx$. All gradient
assumptions refer to distributional derivatives.

\subsection{Local Poincar\'e estimates}

We recall the Poincar\'e estimates in strong and weak Lebesgue spaces.

\begin{lemma}
	\label{lem:poincare}
	For $1<r<\infty$, the following inequalities hold for every
	cube $Q\subset\R^d$ and every $f\in W^{1,1}(Q)$:
	\begin{align}
		\|f-f_Q\|_{L^r(Q)}
		&\lesssim_d \ell(Q)\|\nabla f\|_{L^r(Q)},\label{eq:poincare-strong}\\
		\|f-f_Q\|_{L^r_{\rm weak}(Q)}
		&\lesssim_{d,r} \ell(Q)\|\nabla f\|_{L^r_{\rm weak}(Q)}.
		\label{eq:poincare-weak}
	\end{align}
	If $r>\max\{1,2d/(d+2)\}$, then also
	\begin{equation}\label{eq:poincare-local2}
		\|f-f_Q\|_{L^2(Q)}
		\lesssim_{d,r} \ell(Q)|Q|^{1/2-1/r}
		\|\nabla f\|_{L^r_{\rm weak}(Q)}.
	\end{equation}
	The same assertions hold for balls, with their radius replacing $\ell(Q)$.
\end{lemma}

Inequalities \eqref{eq:poincare-strong} and \eqref{eq:poincare-weak} include the conclusion that $f-f_Q$ belongs to $L^r(Q)$ and $L^r_{\rm weak}(Q)$, respectively, provided $\nabla f$ belongs to the corresponding space.

\begin{proof}
	For the sake of completeness, we first recall the proof of the following pointwise estimate at every Lebesgue point $x\in Q^\circ$:
	\begin{equation}\label{eq:poincare-potential}
		|f(x)-f_Q|\lesssim_d \int_Q |x-y|^{1-d}|\nabla f(y)|\,dy.
	\end{equation}
	For smooth $f$, we integrate the fundamental theorem of calculus along
	the segments from $x$ to $y\in Q$. Convexity gives
	\[
	|f(x)-f_Q|
	\leq\frac1{|Q|}\int_Q\int_0^1
	|x-y|\,|\nabla f(x+t(y-x))|\,dt\,dy.
	\]
	For fixed $t$ we set $z:=x+t(y-x)$. The Jacobian is $t^{-d}$ and
	$|x-y|=|x-z|/t$. If $z$ is in the resulting integration domain,
	then $t\geq |x-z|/\operatorname{diam}(Q)$. Consequently its coefficient
	is at most
	\[
	\frac{|x-z|}{|Q|}
	\int_{|x-z|/\operatorname{diam}(Q)}^1t^{-d-1}\,dt
	\lesssim_d |x-z|^{1-d}.
	\]
	This proves \eqref{eq:poincare-potential}.
	For $W^{1,1}(Q)$ functions, smooth approximation on the cube gives
	the same conclusion almost everywhere: the potentials of the gradient
	errors tend to zero in $L^1(Q)$, so a subsequence converges almost everywhere.
	To obtain the estimate at any Lebesgue point $x\in Q^\circ$, we average this
	almost-everywhere inequality over $B_\varepsilon(x)\subset Q$ and let
	$\varepsilon\downarrow0$. The left side tends to $|f(x)-f_Q|$.
	For $y\ne x$, the averaged kernel satisfies
	\[
	\frac1{|B_\varepsilon(x)|}\int_{B_\varepsilon(x)}|z-y|^{1-d}\,dz
	\lesssim_d |x-y|^{1-d}
	\]
	and tends to $|x-y|^{1-d}$. Thus dominated convergence applies whenever
	the potential on the right of \eqref{eq:poincare-potential} is finite;
	otherwise the desired inequality is trivial.
	
	We extend $|\nabla f|\mathbf1_Q$ by zero and put
	\[
	k_Q(z):=|z|^{1-d}\mathbf1_{\{|z|\leq\operatorname{diam}(Q)\}}.
	\]
	Since $\|k_Q\|_1\lesssim_d \ell(Q)$, Young's inequality proves
	\eqref{eq:poincare-strong}. To obtain its weak counterpart without
	assuming a strong $L^r$ bound, we use the equivalent norm
	\[
	\|F\|_{L^r_{\rm weak}}^{\#}
	:=\sup_{0<|E|<\infty}|E|^{-1+1/r}\int_E|F|\,dx,
	\qquad
	\|F\|_{L^r_{\rm weak}}\leq\|F\|_{L^r_{\rm weak}}^{\#}
	\leq\frac r{r-1}\|F\|_{L^r_{\rm weak}}.
	\]
	The last inequality follows by integrating
	$F^*(t)\leq \|F\|_{L^r_{\rm weak}}t^{-1/r}$, where $F^*$ is the
	decreasing rearrangement of $|F|$; the first follows by testing level sets.
	Fubini and translation invariance give
	\[
	\|k*F\|_{L^r_{\rm weak}}^{\#}\leq\|k\|_1\|F\|_{L^r_{\rm weak}}^{\#}.
	\]
	We apply this to \eqref{eq:poincare-potential} to obtain
\eqref{eq:poincare-weak}. 

For \eqref{eq:poincare-local2}, we choose
$\max\{1,2d/(d+2)\}<t<\min\{r,2\}$. The usual Sobolev--Poincar\'e
inequality and the inclusion of $L^r_{\rm weak}(Q)$ in $L^t(Q)$ give
\[
\|f-f_Q\|_{L^2(Q)}
\lesssim_{d,t}\ell(Q)|Q|^{1/2-1/t}\|\nabla f\|_{L^t(Q)}
\lesssim_{d,r}\ell(Q)|Q|^{1/2-1/r}\|\nabla f\|_{L^r_{\rm weak}(Q)}.
\]
The same argument applies to balls.
\end{proof}

The corresponding estimates also hold in the Birman--Solomyak spaces.

\begin{lemma}\label{lem:poincare-bs}
For $1<r<\infty$ and $E$ equal to either $\ell^r(L^2)$ or
$\ell^r_{\rm weak}(L^2)$, the following inequality holds for every
cube $Q\subset\R^d$ and every $f\in W^{1,1}(Q)$:
\begin{equation}\label{eq:poincare-bs}
\|(f-f_Q)\mathbf1_Q\|_E
\lesssim_{d,r}\ell(Q)\|\mathbf1_Q\nabla f\|_E.
\end{equation}
\end{lemma}

Inequality \eqref{eq:poincare-bs} includes the conclusion that $(f-f_Q)\mathbf1_Q$ belongs to $E$ provided $\mathbf1_Q\nabla f\in E$.

\begin{proof}
	Translations act boundedly on $E$, uniformly in the translation.
	Indeed, after an integer translation, each translated unit cube overlaps
	at most $2^d$ original unit cubes. The \(L^2\) norm of the function on the translated cube is bounded by a sum of at most $2^d$ corresponding local $L^2$ norms. Taking the sequence
	norm proves the assertion. In the weak case we use the equivalent Banach norm
	\[
	\|(a_k)\|_{\ell^r_{\rm weak}}^{\#}
	:=\sup_{\substack{F\subset\mathbb Z^d\\0<\#F<\infty}}
	(\#F)^{-1+1/r}\sum_{k\in F}|a_k|,
	\]
	which is equivalent to the stated weak norm by the same rearrangement
	argument as above. Minkowski's inequality, first on each cube and then
	for the sequence norm, now gives
	\begin{equation}\label{eq:bs-young}
		\|k*F\|_E\lesssim_{d,r} \|k\|_1\|F\|_E.
	\end{equation}
	For general integrals this follows by approximation and the Fatou property;
	in particular it does not require norm continuity of translations in a
	weak space. We apply \eqref{eq:bs-young} to
	\eqref{eq:poincare-potential} to obtain \eqref{eq:poincare-bs}.
\end{proof}

%%%%%%%%%%%%%%%%%%%%%%%%%

\subsection{Approximation modulo constants}

\begin{lemma}\label{lem:approximation}
	We assume that $f\in L^1_{\rm loc}(\R^d)$ and its gradient belongs to one of
	the following spaces $E$:
	\[
	\begin{gathered}
	L^r\text{ or }L^r_{\rm weak}
	\quad\text{with }\max\{1,2d/(d+2)\}<r<\infty,\\
	\ell^r(L^2)\text{ or }\ell^r_{\rm weak}(L^2)\quad\text{with }1<r<2.
	\end{gathered}
	\]
	There are $f_j\in C_c^\infty(\R^d)$ and constants $c_j$ such that
	\begin{equation}\label{eq:approx-local}
		f_j+c_j\longrightarrow f\quad\text{in }L^2_{\rm loc},
		\qquad \sup_j\|\nabla f_j\|_E\lesssim_{d,r} \|\nabla f\|_E.
	\end{equation}
	In the two strong spaces $L^r$ and $\ell^r(L^2)$ one can moreover arrange
	\begin{equation}\label{eq:approx-strong}
		\|\nabla f_j-\nabla f\|_E\longrightarrow0.
	\end{equation}
\end{lemma}

We stress that no convergence in the weak norm is claimed in the weak cases.

\begin{proof}
	The preceding lemmas give $f\in L^2_{\rm loc}$. We choose
	$\chi\in C_c^\infty((-2,2)^d)$ with $\chi=1$ on $[-1,1]^d$, and set,
	for $R\geq1$,
	\[
	Q_R:=(-2R,2R)^d,\qquad c_R:=f_{Q_R},\qquad
	f_R(x):=\chi(x/R)(f(x)-c_R).
	\]
	The product rule and the Poincar\'e estimate in $E$ imply
	\[
	\nabla f_R=\chi(x/R)\nabla f
	+R^{-1}(\nabla\chi)(x/R)(f-c_R),
	\qquad \|\nabla f_R\|_E\lesssim_{d,r} \|\nabla f\|_E.
	\]
	Also $f_R+c_R=f$ on $(-R,R)^d$.
	
	We verify the stronger gradient convergence when $E$ is a strong space.
	The first term in the product rule tends to $\nabla f$ in $E$ by
	dominated convergence. For the second, \eqref{eq:poincare-potential}
	reduces the question to
	\begin{equation}\label{eq:cutoff-error}
		R^{-1}\|k_R*|\nabla f|\|_E\longrightarrow0,
		\qquad k_R(z):=|z|^{1-d}\mathbf1_{\{|z|\leq4\sqrt d R\}}.
	\end{equation}
	The operators $F\mapsto R^{-1}k_R*F$ are uniformly bounded on $E$.
	Bounded functions of compact support are dense in either strong space,
	so it suffices to check \eqref{eq:cutoff-error} for such an $F$.
	Its potential is bounded on bounded sets, and outside a fixed ball
	\[
	|k_R*F(x)|\lesssim_{d,F}|x|^{1-d};
	\]
	the convolution vanishes for $|x|>C_FR$. On the annulus
	$2^m\leq |x|<2^{m+1}$, its $E$ norm is therefore
	$\lesssim_{d,r,F}2^{m(1-d+d/r)}$. This estimate holds for the Birman--Solomyak norm as well:
	the number of unit cubes meeting the annulus is $\lesssim_d2^{md}$, and each
	local $L^2$ norm is $\lesssim_{d,F}2^{m(1-d)}$. Thus
	\[
	\|k_R*F\|_E
	\lesssim_{d,r,F}\left(1+\sum_{0\leq m\leq C_F+\log_2 R}
	2^{m(1-d+d/r)}\right)=o(R),
	\]
	since $1-d+d/r<1$ for $r>1$. This proves
	$\nabla f_R\to\nabla f$ in the strong spaces.
	
	Finally, we take $\rho_\varepsilon$ to be a nonnegative smooth approximate identity
	of compact support and replace $f_R$ by $\rho_\varepsilon*f_R$.
	Convolution preserves the uniform gradient bounds, by Young's inequality
	or \eqref{eq:bs-young}. For fixed $R$ it converges to $f_R$ locally
	in $L^2$. In the strong spaces it also converges in the gradient norm:
	translations are norm continuous there, as follows first for smooth
	compactly supported functions and then by density. We choose $R_j\to\infty$
	and then $\varepsilon_j\downarrow0$ to make these errors tend to zero.
	The functions $f_j:=\rho_{\varepsilon_j}*f_{R_j}$ and constants $c_j:=c_{R_j}$
	have the required properties.
\end{proof}

%%%%%%%%%%%%%%%%%%%%%%%%%%%%%

\subsection{Growth and the definition on test functions}

A uniform local gradient bound above the Sobolev threshold controls
the $L^2$ growth of the function on unit cubes.

\begin{lemma}\label{lem:unit-growth}
	For $\max\{1,2d/(d+2)\}<r<\infty$, if $f\in L^1_{\rm loc}$ and
	\[
	G:=\sup_{k\in\mathbb Z^d}\|\nabla f\|_{L^r(Q_k)}<\infty,
	\]
	then $f\in L^2_{\rm loc}$ and
	\begin{equation}\label{eq:unit-growth}
	\|f\|_{L^2(Q_k)}\lesssim_{d,r}
	|f_{Q_0}|+(1+|k|)G.
	\end{equation}
	Every gradient assumption in Lemma~\ref{lem:approximation} implies
	such a bound for a suitable $r$.
\end{lemma}

\begin{proof}
	The local conclusion and the bound for
	$\|f-f_{Q_k}\|_{L^2(Q_k)}$ follow from
	\eqref{eq:poincare-local2} and $L^r\subset L^r_{\rm weak}$.
	If two unit cubes share a face, a cube of side length at most three
	contains both and meets at most a number of unit cubes bounded solely in terms of the dimension.
	Poincar\'e on that larger cube gives
	$|f_{Q_k}-f_{Q_l}|\lesssim_{d,r}G$ for adjacent cubes $Q_k,Q_l$.
	We join $Q_0$ to $Q_k$ by a path of at most $|k|_1$
	neighboring cubes and add these estimates to obtain
	\eqref{eq:unit-growth}.
	
	For the last assertion, if $\nabla f\in L^p$ or
	$\nabla f\in L^p_{\rm weak}$ with
	$\max\{1,2d/(d+2)\}<p<\infty$, we choose
	$\max\{1,2d/(d+2)\}<r<p$. Since each $Q_k$ has measure one,
	the inclusion $L^p_{\rm weak}(Q_k)\subset L^r(Q_k)$ gives
	\[
	\sup_k\|\nabla f\|_{L^r(Q_k)}
	\lesssim_{p,r}\|\nabla f\|_{L^p_{\rm weak}(\R^d)}<\infty.
	\]
	Here we also use $L^p\subset L^p_{\rm weak}$ in the strong case.
	If instead $\nabla f\in\ell^p(L^2)$ or
	$\nabla f\in\ell^p_{\rm weak}(L^2)$ with $1<p<2$, then
	the sequence $(\|\nabla f\|_{L^2(Q_k)})_{k\in\mathbb Z^d}$
	belongs to $\ell^p$ or $\ell^p_{\rm weak}$, respectively,
	and is therefore bounded. We can thus take $r=2$,
	since $2>\max\{1,2d/(d+2)\}$ in every dimension.
\end{proof}

%%%%%%%%%%%%%%%%%%%%%%

\subsection*{AI declaration}
Generative AI was used to explore estimates proposed by the author and to develop an initial proof strategy; it played a substantial role in constructing the counterexamples. The author checked all mathematical arguments, substantially rewrote an initial manuscript and takes full responsibility for the content.

%%%%%%%%%%%%%%%%%%%%%%%%%%%%%%%%%%%%%%
%%%%%%%%%%%%%%%%%%%%%%%%%%%%%%%%%%%%%%

\bibliographystyle{amsalpha}

\begin{thebibliography}{27}
	
	\bibitem{BHV} B. Bekka, P. de la Harpe, and A. Valette, \textit{Kazhdan’s Property (T)}, Cambridge University Press, 2008.
	
	\bibitem{BS} M.~S.~Birman and M.~Z.~Solomyak, \textit{On estimates of singular numbers of integral operators III. Operators on unbounded domains}, Vest. LSU \textbf{24} (1969); Translation: Vest. LSU Math. \textbf{2} (1975), 9--27.
	
	\bibitem{BKS}
	M. Sh. Birman, G. E. Karadzhov, and M. Z. Solomyak, \textit{Boundedness conditions and spectrum estimates for the operators \(b(X)a(D)\) and their analogs}, Advances in Soviet Mathematics \textbf{7} (1991), 85--106.
	
	\bibitem{Bourgain}
	J. Bourgain, \textit{Some remarks on Banach spaces in which martingale difference sequences are unconditional}, Ark. Mat. \textbf{21} (1983), 163--168.
	
	\bibitem{CGPT}
	J.~M.~Conde-Alonso, A.~M.~Gonz\'alez-P\'erez, J.~Parcet, and E.~Tablate,
	\emph{Schur multipliers in Schatten--von Neumann classes},
	Ann. of Math. (2) \textbf{198} (2023), 1229--1260.
	
	\bibitem{Cwikel}
	M.~Cwikel,
	\emph{Weak type estimates for singular values and the number of bound states
		of Schr\"odinger operators},
	Ann. of Math. (2) \textbf{106} (1977), 93--100.
	
	\bibitem{Davies}
	E.~B.~Davies, \emph{Linear Operators and Their Spectra}, Cambridge Studies
	in Advanced Mathematics \textbf{106}, Cambridge University Press, 2007.
	
	\bibitem{Fr}
	R.~L.~Frank,
	\textit{Cwikel’s theorem and the CLR inequality}, J. Spectr. Theory \textbf{4} (2014), no. 1, 1--21.
	
	\bibitem{Fr2}
	R.~L.~Frank, \textit{The Lieb-Thirring inequalities: recent results and open problems}, In: Kechris, A. (ed.) et al., Nine mathematical challenges. An elucidation. Providence, RI: American Mathematical Society (AMS). Proc. Symp. Pure Math. \textbf{104} (2021), 45--86.
	
	\bibitem{FrLaWe}
	R.~L.~Frank, A.~Laptev, and T.~Weidl,
	\textit{Schr\"odinger operators: Eigenvalues and Lieb--Thirring inequalities}, Cambridge Studies in Advanced Mathematics \textbf{200}, Cambridge University Press, 2023.
		
	\bibitem{FSZ}
	R.~L.~Frank, F.~Sukochev, and D.~Zanin,
	\emph{Endpoint Schatten class properties of commutators},
	Adv. Math. \textbf{450} (2024), 109738.
	
	\bibitem{Ha}
	P. Haj\l asz, \textit{Sobolev spaces on an arbitrary metric space}, Potential Anal. \textbf{5} (1996), 403--415.
	
	\bibitem{Hu}
	D.~Hundertmark, \textit{On the number of bound states for Schr\"odinger operators with operator-valued potentials}, Ark. Mat. \textbf{40} (2002), no. 1, 73--87.
	
	\bibitem{HKRV}
	D.~Hundertmark, P.~Kunstmann, T.~Ried, and S.~Vugalter, \textit{Cwikel’s bound reloaded}, Invent. Math. \textbf{231} (2023), no. 1, 111--167.
	
	\bibitem{KW}
	K. Kazaniecki and M. Wojciechowski, \textit{On Bernstein type quantitative estimates for Ornstein non-inequalities}, Rev. Mat. Iberoam. \textbf{40} (2024), 901--912.
	
	\bibitem{LandauPollak}
	H.~J.~Landau and H.~O.~Pollak,
	\textit{Prolate spheroidal wave functions, Fourier analysis and uncertainty---III: The dimension of the space of essentially time- and band-limited signals},
	Bell System Tech. J. \textbf{41} (1962), 1295--1336.
	
	\bibitem{LSZ}
	G. Levitina, F. Sukochev and D. Zanin, \textit{Cwikel estimates revisited}, Proc. London Math. Soc. (3) \textbf{120} (2020), 265--304.
	
	\bibitem{Li}
	E. H. Lieb, \textit{Bounds on the eigenvalues of the Laplace and Schr\"odinger operators}, Bull. Amer. Math. Soc. \textbf{82} (1976), no. 5, 751--753.
	
	\bibitem{LT}
	E. H. Lieb, W. E. Thirring, \textit{Inequalities for the moments of the eigenvalues of the Schr\"odinger Hamiltonian and their relation to Sobolev inequalities}, Stud. math. Phys., Essays Honor Valentine Bargmann, 269--303 (1976).
	
	\bibitem{Peller80}
	V.~V.~Peller,
	\textit{Hankel operators of class $\mathfrak S_p$ and their applications (rational approximation, Gaussian processes, the problem of majorizing operators)},
	Math. USSR-Sb. \textbf{41} (1982), 443--479; Russian original:
	Mat. Sb. (N.S.) \textbf{113(155)} (1980), 538--581.
	
	\bibitem{Peller24}
	V.~V.~Peller, \textit{Besov spaces in operator theory},
	Russian Math. Surveys \textbf{79} (2024), 1--52.
	
	\bibitem{Ro}
	G. V. Rozenbljum, \textit{Distribution of the discrete spectrum of singular diﬀerential operators}. Dokl. Akad. Nauk SSSR \textbf{202} (1972), 1012--1015. English translation in Soviet Math. Dokl. \textbf{13} (1972), 245--249.
	
	\bibitem{SeSi}
	E. Seiler and B. Simon, \textit{Bounds in the Yukawa\(_2\) quantum field theory: Upper bound on the pressure, Hamiltonian bound and linear lower bound}, Comm. Math. Phys. \textbf{45} (1975), 99--114.
	
	\bibitem{Simon}
	B.~Simon, \emph{Trace Ideals and Their Applications}, second edition,
	Mathematical Surveys and Monographs \textbf{120},
	American Mathematical Society, Providence, RI, 2005.
	
	\bibitem{Simon76}
	B.~Simon, \emph{Analysis with weak trace ideals and the number of bound
	states of Schr\"odinger operators}, Trans. Amer. Math. Soc.
	\textbf{224} (1976), no.~2, 367--380.
	
	\bibitem{Simon79} 
	B. Simon, \textit{Phase space analysis of simple scattering systems: Extensions of some work of Enss}, Duke Math. J. \textbf{46} (1979), no. 1, 119--168.
	
	\bibitem{We}
	T. Weidl, \textit{Cwikel type estimates in non-power ideals}, Math. Nachr. \textbf{176} (1995), 315--334.
\end{thebibliography}

%%%%%%%%%%%%%%%%%%%%%%%%%%%%%%%%%%%%%%
%%%%%%%%%%%%%%%%%%%%%%%%%%%%%%%%%%%%%%

\end{document}